\documentstyle[11pt]{article}
\begin{document}
\input amssym.tex
\def\hgt{\mbox{ht}}
\def\ad{{\mathrm{ad}}\thinspace}
\newcommand{\adj}{\mbox{\em ad\thinspace}}

\def\l{\left}
\def\r{\right}
\def\la{\langle}
\def\ra{\rangle}

\def\hgt{\mbox{ht\thinspace}}
\def\wt{\mbox{wt\thinspace}}

\newtheorem{theorem}{Theorem}[section]
\newtheorem{corollary}{Corollary}[section]
\newtheorem{lemma}{Lemma}[section]
\newtheorem{proposition}{Proposition}[section]
\author{Elizabeth Jurisich}
\date{}
\title{Generalized Kac-Moody Lie algebras, free Lie algebras and
the structure of the Monster Lie algebra}  
\maketitle
\begin{center}
Department of Mathematics, Rutgers University\\ 
New Brunswick, NJ 08903
\end{center}

\section{Introduction}

Generalized Kac-Moody algebras, called Borcherds algebras in
\cite{GT}, were investigated by R. Borcherds in \cite{Bor1}. We show
that any generalized 
Kac-Moody algebra $\frak g$ that has no mutually orthogonal 
imaginary simple roots can be written as $\frak g =\frak u^+ \oplus
(\frak g_J +\frak h) \oplus \frak u^-$, where $\frak u^+$ and $\frak u^-$
are subalgebras  
isomorphic to free Lie algebras with given generators, and $\frak g_J$ 
is a Kac-Moody algebra defined from a symmetrizable Cartan matrix 
(see Theorem \ref{thm:free}). 
There is a formula due to Witt that computes the graded dimension of
a free Lie algebra where all of the generators have been assigned
degree one. It is known that Witt's formula can be extended to other
gradings (e.g., \cite{Bou}). We present a further generalization of the
formula appearing in \cite{Bou}. The denominator identity for
$\frak g$ is obtained
by using this generalization of Witt's formula 
and the denominator identity known for the Kac-Moody 
algebra $\frak g_J$. In this work, we are taking $\frak g$ to be the
algebra defined by the appropriate generators and relations, rather
than the quotient of this algebra by its radical. In particular, our
main result and consequent proof of the 
denominator identity give a new proof that the radical of a generalized
Kac-Moody algebra of the above type is zero. (We use the
fact that the radical of $\frak g_J$ is zero, which is Serre's theorem
in the case that $\frak g_J$ is finite-dimensional; this is the main
case for us.)

The most important application of our work is to the Monster Lie
algebra $\frak m$, defined by R. Borcherds \cite{Bor3}. In fact, we
show that $\frak m = \frak u^+ \oplus \frak g 
\frak l_2 \oplus \frak u^-$, with $\frak u^\pm$ free Lie algebras. 
This result is obtained by applying the above results to a 
generalized Kac-Moody algebra  
$\frak g(M)$ defined from a particular matrix $M$, given by the inner
products of the simple roots of $\frak m$. Theorem
\ref{thm:free} applied to this Lie algebra 
establishes that the subalgebras $\frak n^\pm \subset \frak g(M)$ are
each the semidirect product of a one-%
dimensional Lie algebra and a free Lie algebra on countably many
generators. The Lie algebra 
$\frak g(M)$ is shown to be a central 
extension of the Monster  
Lie algebra $\frak m$  (Theorem \ref{thm:mons}) constructed by R.
Borcherds in  \cite{Bor3}. By Theorem \ref{thm:mons}, the subalgebras
$\frak m^\pm$ in 
$\frak m=\frak m^+\oplus\frak h \oplus \frak m^-$
are isomorphic to the subalgebras $\frak n^\pm \subset \frak g(M)$. In
this way we 
show $\frak m$ contains two large subalgebras $\frak u^\pm$ 
which are isomorphic to
free Lie algebras, and  $\frak m = \frak u^+ \oplus \frak g 
\frak l_2 \oplus \frak u^-$. The denominator identity for $\frak m$ 
(see \cite{Bor3}) is obtained in this paper in the manner described
above for  
more general $\frak g$. In this case $\frak g_J =\frak s \frak l_2$,
and our results 
give a new proof that the central extension $\frak g (M)$ of $\frak m$
has zero radical.

The Monster Lie algebra $\frak m$ is of great interest because R.
Borcherds defines and uses this Lie algebra, along with its denominator
identity, to solve the following problem (\cite{Bor3}):
It was  conjectured by Conway and Norton in \cite{CN} that there
should be 
an infinite-dimensional representation of the Monster simple group
such that the McKay-Thompson series of the elements of the Monster group
(that is, the graded traces of the elements of the Monster group as
they act on the module) are equal to some known modular functions
given in \cite{CN}. After the
``moonshine module'' $V^\natural$ for the Monster simple group was
constructed \cite{FLM1} and many of its properties, including the
determination of some of the McKay-Thompson series, were proven in
\cite{FLM2}, the nontrivial problem of computing the rest of
the McKay-Thompson series of Monster group elements acting  
on $V^\natural$ remained. Borcherds has shown in \cite{Bor3}
that the McKay-Thompson series are the expected modular functions.

In this paper, in preparation for our main result, we include a
detailed treatment of some of Borcherds' work on generalized Kac-Moody
algebras, and of that part of \cite{Bor3} which shows that the
Monster Lie algebra has the properties that we need to prove
Theorem~\ref{thm:mons}. We now explain this exposition.

Some results such as character formulas and a denominator identity
known for Kac-Moody algebras are stated in
\cite{Bor1} for generalized Kac-Moody algebras. We found it necessary
to do some extra work in order to understand fully the precise
definitions and also the reasoning which are implicit in Borcherds'
work on this subject.  V. Kac in \cite{Kac} gives an outline (without
detail) of how to rigorously develop the theory of generalized
Kac-Moody algebras by indicating that one should follow the arguments
presented there for Kac-Moody algebras (see also \cite{HMY}).
Included in \cite{Jur} is a 
detailed exposition of the theory of generalized Kac-Moody algebras,
along the lines of \cite{L}, where the homology results of \cite{GL}
(not covered in \cite{Kac}) are extended to these new Lie algebras.
That this can 
be done is mentioned and used in \cite{Bor3}. The homology result
gives another proof of the character and denominator formulas
(\cite{Jur}). 

We find it appropriate to work with
the extended Lie algebra as in \cite{GL} and \cite{L} (that is, the Lie
algebra with suitable degree derivations adjoined). Alternatively, one
can generalize the theorems in \cite{Kac}. In either of these
approaches the Cartan subalgebra is sufficiently enlarged to make the
simple roots linearly independent and have multiplicity one, just as in
the case of Kac-Moody algebras. Without working in
the extended Lie algebra it does not seem possible to prove the
denominator and character formulas for all generalized Kac-Moody
algebras. This is because the matrix from which we
define the Lie algebra can have linearly dependent columns (as in the
case of $\widehat {\frak s \frak l_2}$); we may even have infinitely
many columns equal. Naturally, when it makes sense to do so, we may
specialize formulas obtained involving the root lattice. In this way
we obtain Borcherds' denominator identity for $\frak m$, and show its
relation to our generalization of Witt's formula.

The crucial link between the Monster Lie algebra and a generalized
Kac-Moody algebra defined from a matrix is provided by Theorem
4.1, which is a theorem given by Borcherds in \cite{Bor2}. Versions
of this theorem also appear in \cite{Bor1} and \cite{Bor3}. Since this
theorem can be stated most neatly in terms of a canonical central
extension of a generalized Kac-Moody algebra (as in \cite{Bor2}) we
include a section on this central extension. 
Theorem~\ref{theorem:hom} roughly says that a Lie algebra with an
``almost positive definite bilinear form'', like the Monster Lie
algebra, is the homomorphic image of a canonical central extension of
a generalized Kac-Moody algebra. The way that
Theorem~\ref{theorem:hom} is stated here and in \cite{Bor2} (as
opposed to \cite{Bor3} where condition {\em 4} is not used)
allows us to conclude that the Monster Lie algebra has a central extension
which is a generalized Kac-Moody algebra defined from a matrix.  
We include in this paper a completely elementary proof of Theorem 4.1. 
This proof is simpler than the argument 
in \cite{Bor1} and the proof indicated \cite{Kac}, which require the
construction of a Casimir operator.
Here equation (11), which follows immediately from the
hypotheses of the theorem, is used in place of the Casimir operator.

The Monster Lie algebra is defined (see \cite{Bor3}) from the
vertex algebra which is the  
tensor product of $V^{\natural}$ and a vertex algebra
obtained from a rank two hyperbolic lattice. This construction is 
reviewed in Section 6.2. The infinite-dimensional 
representation $V^{\natural}$ of the 
Monster simple group constructed in \cite{FLM1} can be given the
structure of a vertex operator algebra, as stated in \cite{Bor0} and
proved in \cite{FLM2}.  The theory
of vertex algebras and vertex operator algebras is used in proving 
properties of the Monster Lie algebra, so the definition of vertex
algebra and a short
discussion of the properties of vertex algebras are given in this
paper. The ``no-ghost'' theorem
of string theory is used here, as it is in \cite{Bor3}, to obtain an 
isomorphism between homogeneous subspaces of the 
Monster Lie algebra and the homogeneous subspaces of $V^\natural$.
A reformulation of the proof of the no-ghost theorem
as it given in \cite{Bor3} and
\cite{Tho} is presented in the appendix of this paper. 

This paper is related to the work of S.J. Kang \cite{Kan}, where a root
multiplicity formula for generalized Kac-Moody algebras is proven
by using Lie algebra homology. We recover Kang's result for the class
of Lie algebras studied in this paper. Other related works include
that of K. Harada, M. Miyamato, and H. Yamada \cite{HMY}, who present
an exposition of generalized Kac-Moody algebras along the lines of
\cite{Kac} (their proof of Theorem 4.1, is the proof in
\cite{Bor1} done in complete detail; see above). The recent work of
the physicists R. Gebert and J. Teschner \cite{GT} explores the module
theory for some basic examples of generalized Kac-Moody algebras.

I would like to thank Professors James Lepowsky and Robert L.
Wilson for their guidance and many extremely helpful discussions.

\section{Generalized Kac-Moody algebras}

\subsection{Construction of the algebra associated to a matrix}

In \cite{Bor1} Borcherds defines the generalized Kac-Moody algebra (GKM)
associated to a matrix. Statements given here without proof have been
shown in detail in \cite{Jur}. In addition to \cite{Bor1}, the reader
may also want to refer to \cite{Kac} where an outline is given for 
extending the arguments given there for Kac-Moody algebras.
The Lie algebra denoted $\frak g'(A)$ in
\cite{Kac}, which is defined from an arbitrary matrix $A$, is equal to the
generalized Kac-Moody algebra $\frak g(A)$, defined below, when 
the matrix $A$ satisfies conditions {\bf C1--C3} given below. 

\noindent{\bf Remark}: In \cite{Bor3}
any Lie algebra satisfying conditions {\em 1--3} of Theorem 4.1 is
defined to
be a generalized Kac-Moody algebra. In this paper the term
``generalized Kac-Moody algebra''  will always mean a Lie algebra
defined from
a matrix as in \cite{Kac}. 
The theory presented here,
based on symmetric rather than symmetrizable matrices,
can be easily adapted to the case where the matrix 
is symmetrizable. We use symmetric matrices in order to be consistent
with the work of R. Borcherds and because the symmetric case is
sufficient for the main applications.

We will begin by constructing a
generalized Kac-Moody algebra associated to a matrix. 

All vector spaces are assumed to be over $\Bbb R$.
Let $I$ 
be a set, at most countable, identified with
${\Bbb Z}_+= \{1,2, \ldots\}$ or with $\{1,2, \ldots ,k\}$. Let 
$A = (a_{ij})_{i,j \in I}$ be a
matrix with entries in ${\Bbb R}$, satisfying the following conditions:
\begin{description}
\item[(C1)]  $A$ is symmetric.
\item[(C2)]  If $ i\neq j$ then $a_{ij}~\leq~0 $.
\item[(C3)]  If $a_{ii} > 0$ then ${2a_{ij} \over a_{ii}} \in \Bbb Z $
for all $j  \in I$. 
\end{description}
  
 Let $\frak g_{0}(A)=\frak g_0$ be the Lie algebra with generators
${h_{i}, e_i, f_i}$, where ${i \in I}$, and the following defining
relations:  
 
\begin{description} 
 \item[(R1)] $\left[ h_{i}, h_{j}\right] =0$
 \item[(R2)] $\left[ h_{i}, e_k\right]   -   a_{ik} e_k =0$
 \item[(R3)] $\l[h_{i}, f_k\r]  +   a_{ik} f_k =0$
 \item[(R4)] $\l[e_i , f_j \r] -   \delta_{ij} h_{i} =0$
\end{description}
 for all $i,j,k \in I$. 
Let
  $\frak h = \sum_{i \in I} {\Bbb R} h_i$. Let 
  $ {\frak n}_0^+ $ be the subalgebra generated by the $ \{e_i\}_{ i \in I}$
    and let $ {\frak n}_0^- $ be the subalgebra generated by the 
   $\{f_i\}_{ i \in I}$. 
The following proposition is proven by the usual methods for
Kac-Moody algebras (see \cite{Jur} or \cite{Kac}):
\begin{proposition}
The Lie algebra $\frak g_0$ has triangular decomposition $\frak g_0 =
{\frak n}_0^- \oplus 
{\frak h} \oplus {\frak n}_0^+$. The abelian subalgebra ${\frak h}$ has  
a basis consisting of $\l\{h_{i} \r\}_{i \in I}$, and
  $\frak n^\pm_0 $ is the free Lie algebra generated by the $e_i$ 
(resp. the $f_i$) $i \in   I$. In particular, $\l\{ e_i, f_i,h_{i}\r\}_{i
\in I}$ is a 
  linearly independent set in $\frak g_0$. \label{prop:dec}
\end{proposition} $\square$
     
 For all $i \neq j$ and $a_{ii} >0$ define
\begin{equation}
 d^+ _{ij} = (\ad e_i)^{1 -2a_{ij} / a_{ii} }e_j \  \in 
{\frak g}_0^+ \label{eq:d1}\end{equation}
\begin{equation}
 d^- _{ij} = (\ad f_i)^{ 1 -2a_{ij} / a_{ii} }f_j \  \in
{\frak g}_0^- \label{eq:d2}\end{equation} 
Let $\frak k_0^{\pm} \subset \frak n_0^{\pm}$ be the ideal of 
$\frak n_0^{\pm}$ generated by the elements: 
\begin{eqnarray}
     d^{\pm}_{ij} & \\
    \l[e_i,e_j\r] & \mbox{ if $ a_{ij}=0$ (in the case of 
$\frak k_0^+$)} \label{e}\\ 
    \l[f_i, f_j\r] & \mbox{ if $ a_{ij}=0$ (in the case of $\frak
k_0^-$).} \label{f} 
\end{eqnarray}
Note that if $a_{ii} >0$, then the elements (\ref{e}) and (\ref{f})
are of type (\ref{eq:d1}) and (\ref{eq:d2}). The subalgebra $\frak k_0
= \frak k^{+}_0 \oplus \frak k^{-}_0$ is an ideal of  
$\frak g_0$ (for details adapt the proof of Proposition \ref{prop:dec2}
in the next section). 

\noindent{\bf Definition 1}: The {\it generalized Kac-Moody 
algebra} $\frak g(A)=\frak g$ associated to the matrix $A$ is the
quotient  
of $\frak g_0$ by the ideal $\frak k_0 = \frak k^+_0 \oplus
 \frak k^-_0$.

\noindent{\bf Remark}: In \cite{Jur} and in $\cite{Kac}$,
a generalized Kac-Moody algebra is constructed as 
a quotient of $\frak g_0$ by its radical (that is, the largest graded
ideal having trivial intersection with $\frak h$). Although 
this fact is not used in this paper, the ideal $\frak k_0$ is equal to
the radical of $\frak g_0$. 
Of course, proving that the radical of the 
Lie algebra $\frak g(A)$ is zero not trivial.
It is shown in \cite{Jur}
that the radical of $\frak g(A)$ is zero using results from \cite{GK},
\cite{Kac} and a proposition proven in \cite{Bor1}. 

Let {$\frak n^+ = \frak n^+_0 / \frak k^+_0$} and {$\frak n^- =
\frak n^-_0 / \frak k^-_0$}. Proposition~\ref{prop:dec} implies that
the 
generalized Kac-Moody algebra has triangular decomposition {$\frak g =  
  \frak n^+ \oplus \frak h \oplus \frak n^-$}.
The Lie algebra $\frak g$ is 
 given by the corresponding generators and relations.
The Lie algebra $\frak g(A)$ is a Kac-Moody algebra when the
matrix $A$ is a generalized Cartan matrix. 

Let $\deg e_i = -\deg f_i = (0,\ldots,0,1,0,\ldots)$ where $1$ appears
in the $i^{th}$ position, and let $\deg h_i =(0,\ldots)$. This induces a Lie
algebra grading by $\Bbb Z^I$ on $\frak g$. Degree
derivations $D_i$ are defined by letting $D_i$ act on the degree
$(n_1,n_2, \ldots)$ subspace of $\frak g$ as multiplication by the
scalar $n_i$. Let $\frak d$ 
be the space spanned by the $D_i$. We extend
the Lie algebra $\frak g$ by taking the semidirect product with 
$\frak d$, so $\frak g^e = \frak d \ltimes
\frak g$. Then ${\frak h}^e = \frak d \oplus \frak h$ is  an
abelian subalgebra of $\frak g^e$, which acts via scalar multiplication
on each space
$\frak g( n_1, n_2, \ldots )$.

Let $\alpha _i \in (\frak h^e)^*$ for $ i \in
I$ be defined by the conditions:
$$[h,e_i] = \alpha_i (h)e_i \mbox{ for all } h \in \frak h^e.$$
Note that $\alpha_j (h_{i}) =  a_{ij}$ for all $i,j \in I$.
Because we have adjoined $\frak d$ to $\frak h$, the $\alpha_i$ are
linearly independent. 

For all $\varphi \in (\frak h^e)^*$ define 
$$\frak g^\varphi = \{x \in \frak g | [h,x] = \varphi (h)x \ \forall h
\in \frak h^e\}.$$ 
If $ \varphi ,\psi\in (\frak h^e)^*$  then $[\frak g^ \varphi , \frak
g^ \psi ] \subset \frak g^ {\varphi + \psi}$.
By definition $e_i \in \frak g^ {\alpha_i},\mbox{ and } f_i \in 
\frak g^{-\alpha_i}$ for all $i \in I$. 
If all $n_i \leq 0$, or all $n_i \geq 0$ (only finitely many nonzero),
it can be shown by using the same methods as for Kac-Moody
algebras that:
$$\frak g^{n_1 \alpha_1 +n_2 \alpha_2 +\cdots} = \frak g(n_1,n_2, \ldots)$$
and $\frak g^0 =\frak h$. Therefore, 
\begin{equation}
\frak g = \coprod_{(n_1,n_2,\ldots)\atop n_i \in \Bbb Z} \frak g^{n_1 \alpha_1
+n_2 \alpha_2 +\cdots} \label{eq:grad}
\end{equation}

\noindent{\bf Definition 2}: The {\it roots} of $\frak g$ are the nonzero
elements $\varphi$
of $(\frak h^e)^* $ such that $\frak g^\varphi \neq 0$.
The elements $\alpha_i$ are {\em simple roots}, and $\frak g^\varphi$ is
the {\em root space} of $\varphi\in (\frak h^e)^*$.

Denote by $\Delta$ the set of roots, $\Delta_+$ the set of {\it positive}
roots i.e., the 
non-negative integral linear combinations of $\alpha_i$.
Let $\Delta_- = -\Delta _+ $
be the set of {\it negative} roots. All of the roots are either
positive or negative.

The algebra $\frak g$  has an
automorphism $\eta$ of order $2$ which acts as $-1$ on $\frak h $
and interchanges the elements $e_i$ and $f_i$. By an inductive
argument, as in \cite{Kac} or \cite{Moo}, we can construct a
symmetric invariant bilinear form on $\frak g$ such that $\frak
g^{\varphi}$ and $\frak g^{-\varphi}$ where $\varphi \in \Delta_+$ are
nondegenerately paired; however, 
the restriction of this form to $\frak h$ can be 
degenerate. There is a character formula for standard modules
of $\frak g$ and a denominator identity (see \cite{Bor1}, \cite{Jur}
and \cite{Kac}):
 \begin{eqnarray}
\prod_{\varphi \in \Delta_+} (1 - e^{\varphi})^{\dim \frak
g^{\varphi}}  & =& \sum_{w \in W } (\det w) \sum_{ \gamma \in
\Omega (0) } (-1)^{l(\gamma)} e^{w(\rho + \gamma)-\rho }
\label{eq:denom} \end{eqnarray}  
where $\Omega(0) \subset \Delta_+ \cup \{0\}$ is the set of all $\gamma\in
\Delta_+\cup \{0\}$ such that $\gamma$ is the sum (of length zero or
greater) of mutually orthogonal imaginary simple roots.

\noindent{\bf Remark}: The denominator formula (\ref{eq:denom}) can 
be specialized to the unextended Lie algebra as
long as the resulting specialization is well defined.

\section{A canonical central extension } 
  
It is useful to consider a certain central extension of the generalized
Kac-Moody algebra.
Working with the central extension defined here
(which is the same as in \cite{Bor2}) will simplify the statement and
facilitate the proof of Theorem~\ref{theorem:hom} below. Given a
matrix $A$ satisfying $\bf C1-C3$ 
let $\hat \frak g$ be the  Lie algebra with 
generators $e_i, f_i, h_{ij}$ for $i,j,k,l \in I$ and relations:
 \begin{description}
\item[(R1$'$)]$\l[h_{ij}, h_{kl}\r]=0$
 \item[(R2$'$)]$\l[h_{ij}, e_k\r] - \delta_{i,j} a_{ik} e_k=0 $ and
   $\l[h_{ij}, f_k\r] + \delta_{i,j} a_{ik} f_k=0$
 \item[(R3$'$)] $\l[e_i , f_j \r] - h_{ij}=0$
 \item[(R4$'$)] $ d^{ \pm}_{ij}=0 \mbox{ for all } i \neq j, a_{ii} > 0$
 \item[(R5$'$)] If $a_{ij} =0$ then $[e_i, e_j]=0$ and $ [f_i,f_j]=0.$  
\end{description}
The elements $d_{ij}^{\pm}$ are
defined by (\ref{eq:d1}) and (\ref{eq:d2}).

We will study this Lie algebra by first considering the Lie algebra
$\hat \frak g_0 $ with generators   
${h_{ij}, e_i, f_i}$, where ${i,j \in I}$, and the defining relations 
($\mbox{\bf R1}')-(\mbox{\bf R3}'$).  
  
Let 
  $\hat \frak h = \sum_{i,j \in I} {\Bbb R} h_{ij} $ and let
  $ {\hat \frak n}_0^+$ be the subalgebra generated by the
$e_i$ $ i
  \in I $,   
 $ {\hat \frak n}_0^-$ be the subalgebra generated by the $ f_i$ $ i \in
   I$.
 We shall prove a version of Proposition~\ref{prop:dec}.

\begin{lemma}
The elements $h_{ij}$ are zero unless the ith and jth columns 
of the matrix $A$ are equal.
\end{lemma}
\noindent {\it Proof\/}:
(also see \cite{Bor2}) The lemma follows from the Jacobi 
identity and relations 
($\mbox{\bf R1}'$), ($\mbox{\bf R2}'$). $\square$

\begin{proposition} The Lie algebra $ {\hat \frak g}_0 ={\hat \frak
n_0}^- \oplus {\hat \frak h} \oplus {\hat \frak n}_0^+ $,
and the abelian Lie algebra $\hat \frak h$ 
has a basis consisting of $\l\{h_{ij} \r\}_{i,j \in I}$ such that the
$i^{th}$ and $j^{th}$ columns of $A = (a_{ij})_{i,j \in I}$ are equal. 
The subalgebra ${\hat \frak n}^\pm_0 $ is the free Lie algebra
generated by the $e_i$ (resp. the $f_i$), $i \in 
  I$. The set $\l\{ e_i, f_i\r\}_{i,j \in I}\cup\l\{ h_{ij}\r\}_{i,j \in S}$ 
is linearly independent where $S=\{(i,j) \in I\times I| a_{ki} = a_{kj}
\mbox{ for all } k \in I\}$.
\label{prop:dec2} 
\end{proposition}  

\noindent {\it Proof\/}: As in the classical case, one constructs a
sufficiently large representation of the Lie algebra. Let 
${\frak h}$ be the span of the elements
$\l\{h_{ij}\r\}_{i,j \in I}$.  
Define $\alpha_j \in {\frak h}^*$ as follows : 
  $$\alpha_j (h_{ik}) = \delta_{ik} a_{ij}.$$
  Let $X$ be the free associative algebra on the symbols $\l\{ x_i \r\}_{i \in
  I}$. Let $ \lambda \in {\frak h}^* $ be such that $ \lambda
(h_{ij})=0$ if 
  $a_{li} \neq a_{lj}$, i.e., unless the $i^{th}$ and $j^{th}$ columns
of $A$ are equal.
  
  We define a representation of the free Lie algebra $\frak g_F$ with
generators
$e_i, f_i, h_{ij}$ on $X$
by the following actions of the generators:
  \begin{enumerate} 
     \item $ h\cdot 1 = \lambda (h)$ $\mbox{ for all } h \in \frak h$
     \item $ f_i\cdot 1 = x_i $ $\mbox{ for all } i \in I$
     \item  $ e_i\cdot 1 = 0 $ $\mbox{ for all } i \in I$
     \item $ h\cdot x_{i_1} \cdots x_{i_r} =( \lambda - \alpha_{i_1} - \cdots
        - \alpha_{i_r}) (h) x_{i_1} \cdots x_{i_r } $ for $h \in \frak h$
     \item $ f_i\cdot x_{i_1} x_{i_2} \cdots x_{i_r} = x_i x_{i_1} x_{i_2}
          \cdots x_{i_r} $
     \item $ e_i\cdot x_{i_1} \cdots x_{i_r} = x_{i_1} e_i \cdot
          x_{i_2} \cdots 
       x_{i_r} + ( \lambda - \alpha_{i_2} - \cdots -\alpha _{i_r}) (h_{ii_1}) 
       x_{i_2} x_{i_3} \cdots x_{i_r} $.  
\end{enumerate}

Let $ {\hat \frak s}$ be the ideal generated by the elements
$\l[h_{ij}, h_{kl}\r]$, $\l[h_{ij}, e_k\r] - \delta_{i,j} a_{ik} e_k$,
$\l[h_{ij}, f_k\r] + \delta_{i,j} a_{ik} f_k$, and $\l[e_i , f_j \r] -
h_{ij}$.  
Now we will show that the ideal $ {\hat \frak s}$ annihilates the module $X$: 
 It is clear that $[h, h'] $ is 0 on $X$ for $h,h' \in {\frak h}$. The
element $[ e_i,f_i] - h_{ij}$ also acts as $0$ on $X$ by the following
computation :
\begin{eqnarray*} 
  [e_i,f_j]\cdot x_{i_1} \cdots x_{i_r} 
     &=& e_i f_j\cdot x_{i_1} \cdots x_{i_r} - f_j e_i\cdot x_{i_1}
\cdots x_{i_r} \\
     &=& e_i\cdot x_j x_{i_1} \cdots x_{i_r} - x_j e_i\cdot x_{i_1} \cdots
     x_{i_r} \\
     &=&(\lambda - \alpha _{i_1} - \cdots -\alpha_{i_r}) (h_{ij}) x_{i_1}
     \cdots x_{i_r} \\
     &=& h_{ij}\cdot x_{i_1} \cdots x_{i_r}. 
\end{eqnarray*}
  By a similar computation $[h_{ij},f_k] + \delta_{ij} a_{ik} f_k$
  annihilates $X$.
  Now consider the action of $[h_{ij},e_k] - \delta_{ij} a_{ik} e_k$ on $X$:
\begin{eqnarray*}  
 [h_{ij}, e_k]\cdot 1 & =& h_{ij} e_k \cdot 1 - e_k h_{ij} \cdot 
  1 \\
      &=& e_k \cdot \lambda(h_{ij})1 \\
      &=&0
\end{eqnarray*}
 and $$\delta_{ij} a_{ik} e_k \cdot 1= 0.$$ 
 Thus $[h_{ij},e_k] - \delta_{ij} a_{ik} e_k$  annihilates $1$.
 Furthermore, $[h_{ij}, e_k] - \delta_{ij} a_{ik} e_k$ commutes with the
 action of $f_l$ for all $l$, as the following computation from
\cite{Bor2} shows:
\begin{eqnarray*} 
\lefteqn {[\ [ h_{ij},e_k ] - \delta_{ij} a_{ik} e_k, f_l \ ]} \\
   & & = [\ [ h_{ij},e_k ], f_l \ ] - [ \delta_{ij} a_{ik} e_k, f_l ]
    = [h_{ij} , h_{kl}] + \delta_{ij} (a_{il} - a_{ik}) h_{kl}.
\end{eqnarray*} 
 By the assumption on $\lambda$, $h_{kl}$ is $0$ on $X$ unless $a_{il}
= a_{ik}$, so that the above is zero. 
 Since any $x_{i_1} \cdots x_{i_r} = f_{i_1} \cdots f_{i_r} \cdot 1$, this 
means that $[h_{ij}, e_k] - \delta_{ij} a_{ik} e_k$ annihilates $X$.
 Now $X$ can be regarded as a $\hat \frak g_0$-module.
The remainder of the proof follows the classical argument. $\square$
 
The following proposition will be
used in the proof of Theorem~\ref{theorem:hom}.

\begin{proposition} In $\hat \frak g_0$, for all $i,j,k \in I$ with $i \neq j$ 
and $a_{ii} > 0$
 $$[e_k,d^-_{ij}] = 0$$ 
 and  $$[f_k,d^+_{ij}] = 0.$$ \label{prop:dij}
\end{proposition}
  
\noindent {\it Proof\/}: It is enough to show the first 
formula.
 
\noindent {\it Case 1\/}: Assume $k \neq i$ and $k \neq j$.
  Since $h_{ki}$ is central if $k \neq i$ 
 \begin{eqnarray}
(\ad e_k) (\ad f_i)x &=& [h_{ki} , x] + [f_i, [e_k, x]] \nonumber \\
             & =& (\ad f_i) (\ad e_k) x \label{eq:ad}
\end{eqnarray}
 so
\begin{eqnarray*} 
[e_k, (\ad f_i)^{{-2a_{ij} / a_{ii}}+1} f_j]& =& 
                (\ad f_i)^{{-2a_{ij} / a_{ii} } +1} [e_k , f_j ] \\
  &=& (\ad f_i)^{{-2 a_{ij} / a_{ii}} +1} h_{kj} = 0 .
 \end{eqnarray*}
The last equality
 holds because $k \neq j$ means $h_{kj}$ is central.
 
\noindent {\it Case 2\/}: Assume $k=i$.
  By assumption $a_{ii}>0$, thus ${e_i, f_i, h_{ii}}$ generate a
Lie algebra isomorphic to $\frak s\frak l_2$. Consider the $\frak s
\frak l_2$-module generated by the weight vector $f_j$. Then if
$a_{ij}=0$ the result follows from the Jacobi identity and the fact
that $[e_i, f_j]= h_{ij}$ is in the center of the Lie algebra. If
$a_{ij}\neq 0$ then 
\begin{eqnarray*}
\lefteqn{\ad e_i(\ad f_i)^{-2a_{ij}/a_{ii} +1}f_j} \\
 &= & (a_{ii}/2) (2a_{ij}/a_{ii} -1)  (-{2a_{ij}/ a_{ii}})
     (\ad f_i)^{-2a_{ij}/a_{ii}}f_j\\
  & & + (-2a_{ij}/a_{ii} +1) (\ad f_i)^{-2a_{ij}/a_{ii}}
(\ad h_{ii})f_j\\
&= &0.
\end{eqnarray*}
 
\noindent {\it Case 3\/}: Assume $k=j$.
 By (\ref{eq:ad})
 \begin{equation}
[e_j, (\ad f_i)^{{-2a_{ij} / a_{ii} } +1 } f_j]=
         (\ad f_i)^{{-2 a_{ij} / a_{ii}}+1} h_{jj}. \label{eq:h}
\end{equation}
Since $[f_i, h_{jj}] = a_{ji} f_i$ it follows immediately that
(\ref{eq:h}) equals zero if $a_{ij} \leq 0$. $\square$
 
  Let ${\frak k}^{ \pm}_0$ be the ideal of ${\hat \frak g_0}^{\pm}$,
respectively, generated by the elements which give relations 
$\mbox{\bf R4}'$ and  $\mbox{\bf R5}'$:
 $$ d^{ \pm}_{ij} \mbox{ for all } i \neq j, a_{ii} > 0$$ 
 $$ [e_i, e_j] \ \mbox{ if }a_{ij} =0 \ \mbox{ for }{\frak k}^+ _0$$
 $$ [f_i,f_j]\ \mbox{ if }a_{ij}=0 \ \mbox{ for }{\frak k}^- _0.  $$
 
\begin{proposition} Define $\frak k_0 = \frak k^+_0 \oplus
\frak k^-_0$.
Then ${\frak k^{\pm}_0}$  
and $\frak k_0$ are ideals of $\hat \frak g_0$.
\end{proposition}
 
\noindent {\it Proof\/}: Let $\hat \frak g_0$ act on itself by the adjoint
representation. 
To see that $\frak k^-$ is an ideal, first consider the action on the
generators $d^-_{ij}$. By Proposition~\ref{prop:dec2} we have
 $$ \sum_{i \neq j} {\cal U}({\hat \frak g}_0) \cdot d^-_{ij} =
\sum_{i \neq j}  
{\cal U}({\hat \frak g}^-_0) {\cal U}({\hat \frak h})
 {\cal U}({\hat \frak g}^+_0) \cdot d^-_{ij}$$
 $$ =\sum_{i \neq j} {\cal U}({\hat \frak g}^-_0) \cdot d^-_{ij} 
 \frak  \subset\frak k^-_0.$$
 The equality holds by Proposition~\ref{prop:dij} and the fact that 
$h(\ad f_i)^N 
 f_j = \lambda (\ad f_i)^N f_j$, where $\lambda$ is a scalar.
 
 We must also consider the action on the generators of the form 
$[f_i,f_j], i \neq j, a_{ij}=0$.
 In this case
\begin{eqnarray*}
[e_k [f_i,f_j]] &=& [h_{ki}, f_j]+[f_i, h_{kj}] \\
     &=&- \delta _{ki} a_{ij} f_j + \delta_{kj} a_{ji} f_i \\
     &=&0
\end{eqnarray*}
i.e. $$ e_k \cdot [f_i, f_j] =0.$$   
 Then by the same argument as above
 $$ \sum_{i \neq j, \ a_{ij}=0} {\cal U}({\hat \frak g}_0) \cdot
[f_i,f_j] \subset {\frak k}^-_0.$$
 By a symmetric argument, $\frak k^+_0$ is an ideal of $\hat \frak g_0$, so
$\frak k_0$ is an ideal. $\square$
 
 The Lie algebra $ \hat \frak g_0 /
\frak k_0 $ is equal to $\hat \frak g$, the Lie algebra 
 defined above by generators and relations. 

\noindent{\bf Remark}: The Lie algebra $\hat \frak g$ is called the 
universal generalized Kac-Moody algebra in \cite{Bor2}. 

Let 
$\frak c$ be the ideal of  
$\hat \frak g$ spanned by the $h_{ij}$ 
where $i \neq j$; note that these elements are central. The Lie
algebra $\hat \frak g$ is a central extension of 
the Lie algebra $\frak g$, because there is an obvious homomorphism
from $\hat \frak g$ to $\frak g$ mapping generators to generators with
kernel $\frak c$. So
$1 \rightarrow \frak c \rightarrow {\hat \frak g} {\buildrel - \over
\rightarrow}  \frak g \rightarrow 1$ is exact. The radical of $\hat \frak g$
must be zero because the radical of $\frak g$ is zero. We have shown
following:

\begin{theorem} The {\it generalized Kac-Moody algebra}
$\frak g$ is isomorphic to ${\hat \frak g}/ \frak c$. 
\end{theorem}

The Lie algebra $\hat \frak g$ can be given a ${\Bbb Z}-$gradation
defined by taking $\deg e_i =-\deg f_i =s_i \in \Bbb Z$ where $s_i
=s_j$ if $h_{ij} \neq 0$, and $\deg h_{ij}=0$. The automorphism $\eta$
is well defined 
on $\hat \frak g$. It follows from Proposition~\ref{prop:dec} and
Proposition~\ref{prop:dec2} that
$\frak n^{\pm}_0 = {\hat \frak n}^{\pm}_0$. 
If $\frak n^+ = {\hat \frak g_0}^+ / \frak k_0^+$ and $\frak n^-
= {\hat \frak g_0}^- / \frak k_0^-$ then 
we have the decomposition $\hat \frak g = \frak n^- \oplus \hat \frak
h \oplus \frak n^+$.

Recall that the Lie algebra $\frak g^e$ has an
invariant bilinear form 
$(\cdot,\cdot)$. It is useful to define an invariant bilinear form
$(\cdot,\cdot)_{\hat \frak g}$ on $\hat \frak g$.
 For $a,b \in \hat \frak g$, let $(a,b)_{\hat
\frak g} = (\bar a, \bar b)$. Note that the span of the $ h_{ij}$
where $i \neq j$ is
in the radical of the form on $\hat \frak g$. The form
$(\cdot,\cdot)_{\hat \frak g}$ is symmetric and invariant because the
form on $\frak g$ has these properties. Grade $\hat \frak g$ by
letting $\deg e_i =1=-\deg f_i$, and $\deg h_{ij}=0$, so $\hat \frak g =
\oplus_{n \in {\Bbb Z}} \hat \frak g_n$, where $\hat \frak g_n$ is
contained in $\frak n^+$ if $n >0$, and $\frak n^-$ if $n < 0$. The
form $(\cdot, \cdot)_{\hat \frak g}$ is nondegenerate on 
$\hat \frak g_n
 \oplus \hat \frak g_{-n}$, because the map given by the central
extension, $x\mapsto \bar x$, is an isomorphism on $\frak n^{\pm}$, and 
the form defined on $\frak g^e$ is nondegenerate on $\frak g^e_m
\oplus \frak g^e_{-m}$ for $m \in \Bbb Z_+$.

\section{Another characterization of GKM algebras}

Theorem \ref{theorem:hom} below is a version of Theorem 3.1 appearing in
\cite{Bor1}. 
Much of the proof of the
following theorem is different than the proof appearing in
\cite{Bor1}. In particular, there is no need to define a Casimir
operator in order to show that the elements $a_{ij}\leq 0$; as seen
below, this follows immediately from condition {\em 3}.

\noindent{\bf Remark:}  In
\cite{Bor3} Borcherds states, as a converse to the theorem, that the
canonical central extensions $\hat \frak g$
satisfy conditions {\em 1--3} below, although we note that the canonical
central extension of a generalized Kac-Moody algebra does not have
to satisfy condition {\em 1} for some matrices. For example, if we
start with the infinite matrix whose 
entries are all $-2$, then all of the $e_i$ must have the same degree
because of condition {\em 3}. Therefore,
there is no way to define a $\Bbb Z-$grading of $\hat \frak g$ so that
$\hat \frak g_i$ is both finite-dimensional and satisfies condition
{\em 3}.

We also note that the kernel of the map $\pi$ appearing in the
following theorem can be strictly larger that the span of the
$h_{ij}$; cf. the statement of Theorem 4.1 in \cite{Bor3}. 

\begin{theorem}[Borcherds] Let $\frak g$ be a Lie algebra satisfying the
following conditions:
\begin{enumerate}
 \item  $\frak g$ can be $\Bbb Z$-graded as $\coprod_{i\in {\Bbb
Z}} \frak g_i$, $\frak g_i$ is finite dimensional if $i \neq 0$, and
$\frak g$ is diagonalizable with respect to $\frak g_0$.
\item  $\frak g$ has an involution $\omega$ which maps $\frak g_i$
onto $\frak g_{-i}$ and acts as $-1$ on noncentral elements of 
$\frak g_0$, in particular $\frak g_0$ is abelian.
\item  $\frak g$ has a Lie algebra-invariant bilinear form $(\cdot,\cdot)$,
invariant under $\omega$, such that $\frak g_i$ and $\frak g_j$ are
orthogonal if $i \neq -j$, and such that the form $(\cdot,\cdot)_0$,
defined by $(x,y)_0 = -(x, \omega (y))$ for $x,y \in \frak g$, is positive
definite on $\frak g_m $ if $m \neq 0$.
\item  $\frak g_0 \subset [\frak g,\frak g]$. 
\end{enumerate} 
Then there is a central extension $\hat {\frak g}$ of a generalized
Kac-Moody algebra and a homomorphism, $\pi$, from $\hat
{\frak g}$ onto 
$\frak g$, such that the kernel of $\pi$ is in the center of $\hat
{\frak g}$. 
\label{theorem:hom}
\end{theorem}

\noindent {\it Proof\/:}
Generators of the Lie algebra $\frak g$ 
are constructed, as in \cite{Bor1}, as follows: For $m >0$, let $\frak l_m$ 
be the subalgebra of $\frak g$ generated by the $\frak g_n$ for $0
<n<m$, and let $\frak e_m$ be the orthogonal complement of $\frak l_m$
in $\frak g_m$ under $(\cdot,\cdot)_0$. To see that $\frak
e_m$ is invariant under $\frak g_0$ let $x \in \frak g_0,\ y \in \frak
e_m,$ and $ \ z \in \frak l_m$. Then $[x,z] \in 
\frak l_m$, so $ (y,[x,z])_0 =0$. Since the form $(\cdot,\cdot )_0$
satisfies $([x,y],z)_0 = -(y,[\omega (x),z])_0$, i.e., is contravariant,
we have $([x,y],z)_0 =0$, which implies that  $[x,y] \in \frak e_m$.
The operators induced by the action of $\frak g_0$ on $\frak e_m$ commute,
so we can construct a basis of $\frak e_m$ consisting of weight
vectors with respect to $\frak g_0$. The form $(\cdot,\cdot)_0$ is positive
definite on $\frak e_m$ so an orthonormal basis can be constructed.
Contravariance of the form ensures that this orthonormal basis also 
consists of weight vectors. The union of these bases for all the 
$\frak e_m$'s can be indexed by $I = {\Bbb Z_+}$, in any order, and
will be denoted 
$\{e_i\}_{i\in I}$. Each $\frak g_n,\ n>0$, is in the Lie algebra
generated by $\{e_i\}_{i\in I}$, as is seen by the following induction
on the degree, $n$: For $n=1$, $\frak g_1 = 
\frak e_1$. Now assume that the $\frak g_n$ for all $0<n<m$ are
contained in the Lie algebra generated by the $e_i$, $i\in I$. The
finite dimensional space $\frak g_m$ 
decomposes under $(\cdot,\cdot)_0$ as $\frak g_m =\frak e_m \oplus \frak
e_m ^\perp $, where $\frak e_m^\perp = \frak l_m \cap \frak g_m$. 
By the induction assumption, 
$\frak e_m^\perp$ is generated by some of the $e_i$'s, and by
construction $\frak e_m$ has a basis consisting of $e_i$'s. 
Define $f_i = -\omega
(e_i)$, and $h_{ij} = [e_i, f_j]$. The $\frak g_n$ where $n <0$
are generated by the $f_i$, $i \in I$. 

The elements $h_{ij}$ can be nonzero only
when $\deg e_i = \deg e_j$. This is because if $\deg e_i > \deg e_j$, which
can be assumed without loss of generality, then 
$[e_j,[e_i,f_j]] \in \frak l_{\deg e_i} $. Thus $([e_i, f_j],
[e_i,f_j])_0 =0$ by contravariance, and $[e_i, f_j]=0$ by the positive
definiteness of $(\cdot, \cdot)_0$. Therefore, all of the $h_{ij}$ are
in $\frak g_0$. By assumption {\em 4}, $\frak g_0$ is generated by
the $h_{ij}$, $i,j \in I$.

Let
$k \in$ rad$(\cdot,\cdot)$ (so $k \in \frak g_0$). 
Then $([k,g],[k,g])_0 = (k,[\omega
(g),[k ,g]])_0 =0$ for all $g \in \frak g_n$, $n \neq 0$.
Thus $[k,g] =0$ by positive
definiteness, and since $k \in \frak g_0$, $k$ must also commute with 
$\frak g_0$. Therefore the radical of the form $( \cdot,\cdot)$ is in the
center of $\frak g$.

If $ i\neq j $ then $[e_i, f_j]=h_{ij}$ is contained 
in the center of $\frak g$ because it is in the radical of the form
$(\cdot, \cdot)$. To see this,
consider $([e_i,f_j], x)$ for a homogeneous $x \in \frak g$. We know
by assumption that 
this is zero if $x \notin \frak g_0$. If $x \in \frak g_0$ then 
$([e_i,f_j],x) =(e_i, [f_j,x])= c(e_i, e_j)_0 =0$ for some real number $c$,
as $f_j$ is a weight vector for $\frak g_0$, and the $e_k$'s are
orthogonal. Thus $([e_i,f_j],x)=0 \mbox{ for all } x \in \frak g$.

Now it must be shown that the generators constructed above satisfy the
the relations  $\mbox{\bf R1}'$-$\mbox{\bf R5}'$
of the central extension of a generalized Kac-Moody algebra, and
the symmetric matrix with entries 
$(h_{ii},h_{jj}) = a_{ij}$ satisfies conditions {\bf C2 - C3}.

That $\mbox{\bf R1}'$ (that is, $[h_{ij}, h_{kl}] = 0$), holds is obvious 
because
$\frak g_0$ is abelian, and $\mbox{\bf R3}'$ is true by definition. 

The relations $\mbox{\bf R2}'$
are proven by the following argument (cf. \cite{Bor1}):
By construction, 
$e_i$ is a weight vector of $\frak g_0$, thus for some real number $c$,
$[h_{lm},e_i] = ce_i$ and
\begin{eqnarray*}
  c=c(e_i,e_i)_0 & = & ([h_{lm},e_i],e_i)_0 \\
                 & =& ([h_{lm}, e_i],f_i) \\
                 & =&  (h_{lm},[e_i,f_i])  = (h_{lm}, h_{ii}) .   
\end{eqnarray*}
Since $h_{ij}$ for $i \neq j$ is in the radical
of the form $(\cdot, \cdot)$, we have $c = \delta_{lm}a_{li}$. Applying
$\omega$ shows the relation for $f_i$. 

To show $\mbox{\bf R5}'$ and {\bf C2} let $i \neq j$.
We must prove $a_{ij} \leq 0$,
and if $a_{ij} =0$ then $[e_i,e_j]=0$ and $[f_i,f_j]=0$. 
Consider the element $[e_i, e_j] \in \frak g$, we have 
 \begin{eqnarray}
 ([e_i,e_j],[e_i,e_j])_0 & = & -(e_i,[f_j,[e_i,e_j]])_0 \nonumber \\
                         & = & -a_{ij}(e_i,e_i)_0.  \label{eq:a}
\end{eqnarray}
The last equality follows from 
\begin{eqnarray*}
 [f_j,[e_i,e_j]] &= & -[e_j,[f_j,e_i]]- [e_i,[e_j,f_j]] \\
                           &=& a_{ij} e_i. 
\end{eqnarray*}
By equation (\ref{eq:a}) and the
positive definiteness of the form $(\cdot,\cdot)_0$ on $\frak g_m ,\ m
\neq 0,$ we have $a_{ij} \leq 0$, and $a_{ij}=0$ if and only if
$[e_i,e_j]=0$. By applying $\omega$ we also show $[f_i,f_j] =0$ in this case.


For $a_{ii} >0$ the Lie algebra generated by ${e_i ,f_i, h_{ii}}$ is
isomorphic to $\frak s \frak l _2$, rescaling $e_i$ and $h_{ii}$ by 
$2 / a_{ii}$. For each $j \in I$ the element $f_j$ generates an $\frak
s \frak l_2$-weight module, the weights must all be integers so
$[h_{ii}, f_j] =(-2a_{ij} / a_{ii}) f_j$ implies that $2a_{ij} / a_{ii} 
\in {\Bbb Z}$, thus {\bf C3} is satisfied. Proposition~\ref{prop:dij}, 
which holds for $\frak g$, shows 
that $[e_k, d_{ij}^-] =0$ where $d_{ij}^- = (\ad f_i)^{n+1} f_j$ and
$n = -2a_{ij}/a_{ii}$ (note that $n$ is positive). 
Contravariance of the form gives us
$$((\ad f_i)^{n+1}f_j ,(\ad f_i)^{n+1}f_j)_0 = 
    ((\ad f_i)^{n}f_j ,[e_i,(\ad f_i)^{n+1}f_j])_0 
     = 0, $$
so that $(\ad f_i)^{1 - 2a_{ij} / a_{ii}}f_j= 0$ by the positive
definiteness of $(\cdot, \cdot)_0$.  Applying $\omega$ to $d_{ij}^-$
gives the relation for $(\ad e_i)^{1 - 2a_{ij} / a_{ii}}e_j$. This 
shows the relations $\mbox{\bf R4}'$.

Denote by $\hat {\frak g}$ the canonical central extension, defined in 
Section $ 3$,
of the  generalized Kac-Moody algebra 
associated to the matrix $ a_{ij}=(h_{ii},h_{jj})$.
Define a homomorphism $\pi : \hat {\frak g} \rightarrow \frak g$,
taking the generators ${e_i, f_i, h_{ij}}$ in $\hat {\frak g}$ to the
generators ${e_i, f_i, h_{ij}}$ in $\frak g$. Since the generators of
$\frak g$ have been shown to satisfy the relations of $\hat \frak g$
the map $\pi$ is a homomorphism from $\hat {\frak g}$ onto 
$\frak g$. The bilinear form, $(\cdot,\cdot)_{\hat {\frak g} }$,
on $\hat {\frak g}$ satisfies
$$(e_i,f_j)_{\hat {\frak g}} = \delta_{ij}= (e_i,f_j)$$
$$(h_{ii},h_{jj})_{\hat {\frak g}} = a_{ij} =(h_{ii},h_{jj}).$$
The $h_{ij}$ with $i \neq j$ are in the radical of $(\cdot, \cdot
)_{\hat \frak g}$. 
Thus $(x,y)_{\hat {\frak g}} = (\pi(x),\pi (y))$ for $x,y \in \hat
\frak g$, because
$(x,y)_{\hat {\frak g}} $ can be reduced using invariance and the
Jacobi identity to some polynomial in the $(e_i,f_j)$ and $(h_{ii},h_{jj})$.

Now we determine the kernel of the map $\pi$. Let $a \in {\hat
{\frak g}}$ such that $a \neq 0$ and $\pi (a) =0$ 
in $\frak g$. Recall the grading $\hat \frak g = \oplus_{n \in \Bbb
Z} \hat \frak g_n$ and the decomposition ${\hat {\frak g}}=  {\frak
n}^+ \oplus \hat {\frak h}  \oplus  {\frak n}^-$, thus can write $a=
a_+ + a_0 +a_-$ 
where $a_{\pm} \in  {\frak n}^{\pm} \ ,a_0 \in \hat {\frak h}$.
Thus $\pi (a) = \pi (a_+) + \pi 
(a_0) +\pi (a_-) =0$, which is still a direct sum in $\frak g$. Therefore, 
$\pi (a_+)=0 \ ,\pi(a_0)=0 \ , \pi(a_-)=0$. Assume that $a_+$ is
homogeneous and nonzero, then for some $n >0$, $a_+ \in \hat {\frak
g}_{n}$ and $(a_+, x)_{\hat {\frak g}}  =
(\pi (a_+),\pi (x))=0 \mbox{ for all } x \in \hat {\frak g}_{-n}$.
Since $(\cdot,\cdot)_{\hat {\frak g}}$ is nondegenerate on $\hat
{\frak g}_{n }\oplus \hat {\frak g}_{-n}$ we have $a_+ =0$.
By a similar argument $a_- =0$. Thus $a =a_0 \in \hat {\frak h}$, and
$[a,h]=0$ for all $h \in \hat \frak h$. Since $\pi(a) =0$, we have
$$(a,h_{ii})_{\hat {\frak g}}=(\pi (a),\pi (h_{ii})) =0 \mbox{ for all }
i$$
so $$[a,e_i]=(a,h_{ii})_{\hat {\frak g}} e_i =0 \mbox{ for all } i$$
similarly $[a,f_{i}] =0$ for all $i \in I$. Thus $a $ is in the center of
$\hat {\frak g}$. $\square$

If the radical of the form $(\cdot,\cdot)$ is zero then the elements
$h_{ij}$ are all zero and we have a homomorphism to $\frak g$ from a
generalized 
Kac-Moody algebra, for which character and denominator formulas
have been 
established. By construction of the monster Lie algebra in
\cite{Bor2}, 
which is also discussed later in this paper, the
radical of the invariant form is zero, so that the following corollary
will apply to this algebra.

\begin{corollary} Let $\frak g$ be a Lie algebra satisfying the
conditions in Theorem~\ref{theorem:hom}. If the radical of the form on
$\frak g$ is zero then there is a generalized Kac-Moody
algebra $\frak l$ such that $\frak l / \frak c = \frak g$, where
$\frak c$ is the center of $\frak l$. \label{cor:bor} 
\end{corollary}

\section{Free subalgebras of GKM algebras}

\subsection{Free Lie algebras}

Denote by $L(X)$ the free Lie algebra on a set $X$, if $W$ is 
a vector space with basis $X$ then we also use the notation 
$L(W) = L(X)$. We will assume that $X$ is
finite or countably infinite.
Let $\nu =(\nu_1,\nu_2,\ldots)$ be an $m-$tuple where $m \in \Bbb Z_+
\cup {\infty}$, and $\nu_i \in \Bbb N$ satisfy
$\nu_i =0$ for $i$ sufficiently large.
We will use the notation $|\nu| =\sum_{i=1}^m \nu_i$.
If $T_i$ are indeterminates, let
$$T^{\nu} = \prod_{i \in \Bbb Z_+}T_i^{\nu_i} \in 
\Bbb R[T_1,\ldots, T_m] \mbox{ or } \Bbb R[T_1,T_2,\ldots] .$$ 

We will consider gradings of the Lie algebra $L(X)$ of the following
type: Let $\Delta = \Bbb Z^m$ and assign to each element 
$x \in X$ a degree $\alpha \in 
\Delta$ i.e., specify a map $\phi : X \rightarrow \Delta$. Let
$n_\alpha$ be the number of elements of $X$ of 
degree $\alpha$ and assume this is finite for all $\alpha \in \Delta$.
This defines a grading of $L(X)$ which we denote $L(X) =
\coprod_{(\alpha \in \Delta)} L^{\alpha}(X)$. 
Let $d(\alpha) = \dim L^{\alpha}(X)$. We assume that $d(\alpha)$ is finite.

An example of this type of
grading is the {\em multigradation} from Bourbaki \cite{Bou}, which is
defined as follows:
Enumerate the set $X$ so that $X 
=\{ x_i\}_{i \in I}$ where $I = \Bbb Z_+$ or $\{1,\ldots,m\}$. Denote
by $\Delta_{|X|}$ the group of $|X|-$tuples.
Let $\sigma:X \rightarrow \Delta_{|X|}$ be defined by $\sigma(x_i)=
\epsilon_i =(0,\ldots,0,1,0,\ldots)$ 
where the $1$ appears in position $i$. The map $\sigma$ defines the
multigradation of $L(X)$.

Given such a grading determined by the map $\phi$ and the group
$\Delta = \Bbb Z^m$, if 
$\Delta' = \Bbb Z^n$ and if $\psi : \Delta 
\rightarrow \Delta'$ is a homomorphism, then $\psi \circ \phi : X \rightarrow
\Delta'$ also defines 
a grading of the Lie algebra $L(X)$. It is clear that if $\alpha \in
\Delta'$ then $d(\alpha)
=\sum_{\beta \in \Delta \atop \psi (\beta) = \alpha} d(\beta)$ as
$L^{\alpha}(X) = 
\coprod_{\beta \in \Delta \atop \psi (\beta) = \alpha} L^{\beta}(X)$.
Notice that $d(\alpha)$ is zero unless $\alpha$ has all nonnegative or
all nonpositive entries, if the degree of each element of $X$ has this 
property.
The map $\psi$ induces a 
homomorphism from $\Bbb R[T_1,\ldots, T_m]$ to $\Bbb R[T_1, \ldots
T_n]$, where if $\alpha \in \Delta$ then $T^\alpha \mapsto T^{\psi(\alpha)}$.

\begin{proposition}
Let $L(X) = \coprod_{(\alpha \in \Delta)} L^{\alpha}(X)$ be a grading
of the above type. Then
\begin{equation}
1- \sum_{\alpha \in \Delta} n_{\alpha}T^{\alpha} = \prod_{\alpha \in
\Delta\backslash\{0\}} (1- T^{\alpha})^{d(\alpha)}. \label{eq:wil}
\end{equation} \label{prop:will}
\end{proposition}

\noindent{\it Proof}:
If we
consider the multigradation of $L(X)$ by $\Delta_{|X|}$, so $L(X) =
\coprod_{\beta \in\Delta_{|X|}}L^{\beta}(X)$, the following formula
is proven in \cite{Bou} 
(for finite $X$, but this immediately implies the result for an
arbitrary $X$): 
\begin{equation}
1-\sum_{i \in I}T_i = \prod_{\beta \in \Delta_{|X|}\backslash \{0\}}(1-
T^{\beta})^{d(\beta)}.\label{eq:bou} \end{equation} 
This implies a formula for the more general type of grading given 
by a map $\phi : X \rightarrow \Delta$, as above. Any such map
satisfies $\phi =\phi' \circ \sigma $ where $ \phi' : \Delta_{|X|}
\rightarrow \Delta$ is given by:
$$\epsilon_i \mapsto \phi(x_i).$$
Applying the homomorphism $\phi'$ to the identity (\ref{eq:bou})
gives the proposition.$\square$

Some of our results  will follow from the elimination theorem in
\cite{Bou}, which is restated here for the convenience of the reader.

\begin{lemma}[elimination theorem] Let $X$ be a set, $S$ a subset of
$X$ and $T$ the set of sequences $(s_1,\ldots,s_n,x)$ with $n \geq 0,
s_1,\ldots,s_n$ in $S$ and $x$ in $X\backslash S$.
\begin{description}
\item[(a)] The Lie algebra $L(X)$ is the direct sum as a vector space
of the subalgebra 
$L(S)$ of $L(X)$ and the ideal $\frak a$ of $L(X)$ generated by
$X\backslash S$.
\item[(b)] There exists a Lie algebra isomorphism $\phi$ of $L(T)$
onto $\frak a$ which maps $(s_1,\ldots,s_n,x)$ to $(\adj s_1\cdots   
\adj s_n) (x)$. \end{description}
\end{lemma}$\square$


As Bourbaki \cite{Bou} does for two particular gradings, we
obtain formulas for computing the 
dimension of the homogeneous subspaces of a free Lie algebra $L(X)$ 
graded as above by a group $\Delta$. The formulas 
derived here relate the dimension of the piece of degree $\alpha$ with
the number of generators in that degree (which is assumed finite).
If $\beta \in \Delta$ can be partitioned $\beta =
\sum a_{\alpha} \alpha$, where $a_{\alpha} \in \Bbb N$ and $\alpha \in
\Delta$, then define the partitions $P(\beta, j) = \{a
=(a_{\alpha})_{\alpha \in \Delta}| \beta =\sum a_{\alpha}\alpha, 
 |a|= j \}$ and $P(\beta) = \cup_j P(\beta, j)$.
Taking $\log$ of both sides of formula (\ref{eq:wil}) leads to the
equations: 
$$-\log (1-\sum_{\alpha \in \Delta}n_{\alpha}T^{\alpha}) = \sum_{j\geq 1}
{1\over j}(\sum_{\alpha \in \Delta} n_{\alpha} T^{\alpha})^j
 = \sum_{j \geq 1} {1 \over j} \sum_{\beta \in
\Delta}\sum_{a \in P(\beta ,j)}{|a|! \over a!}
\prod_{\alpha \in \Delta} n_{\alpha}^{a_{\alpha }} T^{\beta}$$
and 
$$-\sum_{\alpha \in \Delta} d(\alpha)\log(1-T^{\alpha}) =
\sum_{\alpha \in \Delta ,k \geq 1}{1 \over k} d({\alpha})T^{k\alpha}
 = \sum_{\beta} \sum_{k |\beta} {1\over k} d(\beta/k) T^{\beta}.$$
Thus if $\gamma | \beta$ means $k\gamma =  \beta$ then
$$\sum_{\gamma | \beta}{\gamma \over \beta} d(\gamma) = 
\sum_{a \in P(\beta)}{(|a|-1)! \over a!}
 n^a $$
where $ n^a $ denotes the product of the $n_{\alpha}^{a_\alpha}$.
Applying the M\"{o}bius inversion formula gives:
\begin{proposition} Let $\Delta$ be a grading of the free Lie algebra
$L(X)$ as in Proposition~\ref{prop:will}. If $d(\beta)$ is the
dimension of $L^\beta (X)$ for $\beta \in \Delta$ then
\begin{equation}
d(\beta)=\sum_{\gamma|\beta}({\gamma \over \beta})\mu({\beta \over
\gamma})\sum_{a \in P(\gamma )}{(|a|-1)! \over a!} n^a. \label{eq:dim} 
\end{equation} 
\end{proposition}

\subsection{Applications to generalized Kac-Moody algebras}

In this section we will show that certain generalized Kac-Moody
algebras contain large subalgebras which are isomorphic to
free Lie algebras. We will apply the results of the preceding section
to these examples.

We will begin with an easy example. Let $A$ be a matrix
satisfying conditions $\mbox{\bf C1-C3}$ which has no $a_{ii} >0$, and all
$a_{ij} < 0$ for $i \neq j$ (this means that imaginary simple roots
are not mutually orthogonal). In this case, the generalized Kac-Moody
algebra $\frak g(A)$ is equal to $\frak g_0(A)$. By
Proposition~\ref{prop:dec}, $\frak g = \frak n^+ \oplus \frak h \oplus
\frak n^- $ where $ \frak n^{\pm}$ are the free Lie algebras on the
sets $\{e_i\}_{i \in I}$ and $\{f_i\}_{i \in I}$ respectively.

\noindent{\bf Remark}: Formula (\ref{eq:wil}) can be applied to the
root grading 
of $\frak g(A)$ to obtain the denominator identity for $\frak g(A)$.
The root multiplicities are given by (\ref{eq:dim}).

If the imaginary simple roots of a generalized Kac-Moody algebra
are not mutually orthogonal then $\frak n^\pm$ are not in general 
free, but we
will show they contain ideals which are isomorphic to free Lie
algebras. First 
we will set up some notation. Let $J \subset I$ be the set $\{ i \in I
| \alpha_i \in \Delta_R \} = \{ i \in I | a_{ii} >0\}$. Note that the
matrix $(a_{ij})_{i,j \in J}$ is a generalized Cartan matrix, let
$\frak g_J$ be the Kac-Moody algebra associated to this matrix.
Then $\frak g_J = \frak n^+_J \oplus \frak h_J \oplus \frak n_J^-$, and
$\frak g_J$ is isomorphic to the subalgebra of $\frak g(A)$ generated 
by $\{e_i, f_i\}$ with $i \in J$. 

\begin{theorem}
 Let $A$ be a matrix satisfying conditions {\bf C1-C3}. Let $J$ and
$\frak g_J$ be as above. Assume 
that if $i,j \in I \backslash J$ and $i \neq j$ then $a_{ij}<0$. 
Then $$\frak g(A) = \frak u^+ \oplus (\frak g_J + \frak h) \oplus \frak u^-,$$ 
where
$\frak u^- = L(\coprod_{j \in I\backslash J}{\cal U}(\frak n^-_J)\cdot f_j) $
and
$\frak u^+= L(\coprod_{j \in I\backslash J}{\cal U}(\frak n^+_J)\cdot e_j ) $.
The ${\cal U}(\frak n^-_J)\cdot f_j $ for $j \in I\backslash J$ are integrable 
highest weight $\frak g_J$-modules, and the ${\cal U}(\frak n^+_J)\cdot e_j $ 
are integrable lowest weight $\frak g_J$-modules.
\label{thm:free}
\end{theorem}

Note that the conditions on the $a_{ij}$ given in the theorem are
equivalent to the statement that the Lie algebra has no mutually
orthogonal imaginary simple roots. 

\noindent{ \it Proof\/}: We will consider $\frak n^+$; the case of
$\frak n^-$ is shown by a similar argument or by applying the
automorphism $\eta$. 

By the construction in Section 2, we have
  $$ \frak n^+ = L(\{e_{i}\}_{i \in I})/ \frak k_0^+,$$
where  
$\frak k_0^+$ is generated as an ideal of $L(\{e_{i}\}_{i \in I})$ 
by the elements
$$\{ (\ad e_{i})^{1-2a_{ij}/a_{ii}}e_{j}\ | \ i,j \in J, i\neq j \}$$ and
$$\{ (\ad e_{i})^{1-2a_{ij}/a_{ii}}e_{j}\  | \ i \in J, j \in I 
\backslash J\}.$$
This is because there are no elements of type (\ref{e}).
Apply the elimination theorem to the free Lie algebra $L(\{e_{i}\}_{i
\in I})$ with $S =J$. Thus
$$L(\{e_{i}\}_{i \in I})= L(\{e_{i}\}_{i \in J}) \ltimes \frak a$$
where the ideal $\frak a$ is isomorphic to the free Lie algebra on the
set $X =\{ \ad e_{i_1}\ad e_{i_2}\cdots \ad e_{i_k} e_j\ |\  j \in
I \backslash J \mbox{ and } i_m \in J\}$. Let $W$ denote the vector 
space with basis $X$, so that $\frak a \cong L(W)$.
Observe that as an $\frak h^e$-module 
$W \cong \coprod_{j \in I \backslash J} {\cal U}(\frak l) e_j $,
where $\frak l$ denote the free Lie algebra
$L(\{e_i\}_{i \in J})$.

For each fixed $j \in I\backslash J$ consider the submodule
$$R_j = \coprod_{i \in J} {\cal U} (\frak l) 
(\ad e_i)^{1- 2a_{ij}/a_{ii}}e_j \subset {\cal U}(\frak l) e_j .$$ 
Thus, identifying quotient spaces with subspaces of $W$,
\begin{eqnarray*}
W &=& \coprod_{j \in I \backslash J} ( {\cal U}(\frak l) e_j /R_j \oplus R_j)\\
  &=& \coprod_{j \in I \backslash J}{\cal U}(\frak l) e_j/R_j \oplus 
        \coprod_{j \in I \backslash J} R_j
\end{eqnarray*}

Now apply the elimination theorem to the Lie algebra 
$L(X) = L(W)$, choosing a basis of $W$ of the form 
$S_1 \cup S_2$ where $S_1$ is a basis of the vector space
$\coprod_{j \in I \backslash J}{\cal U}(\frak l) e_j/R_j$ and
$S_2$ is a basis of $\coprod_{j \in I \backslash J} R_j$. 
Obtaining

$$L(W) = L(\coprod_{j \in I \backslash J}{\cal U}(\frak l) e_j/R_j ) 
              \ltimes \frak b$$
where $\frak b$ is the ideal of $L(W)$ that is generated by $S_2$,
i.e., by $\coprod_{j \in I \backslash J} R_j$. 
So
$$L(\{e_{i}\}_{i \in I})= L(\{e_{i}\}_{i \in J}) \ltimes  
L(\coprod_{j \in I \backslash J}{\cal U}(\frak l) e_j/R_j ) 
\ltimes \frak b.$$ 

Let $\frak k_J^+$ be the ideal of $\frak l$ generated
by $$\{ (\ad e_{i})^{1-2a_{ij}/a_{ii}}e_{j}\ | \ i,j \in J, i\neq j \},$$
then ${\cal U}( \frak n_J^+) = {\cal U} (\frak l/\frak k_J^+) =
{\cal U} (\frak l)/ {\cal K}$, where $\cal K$ denotes the ideal of 
${\cal U}(\frak l)$ generated by $\frak k_J^+ \subset {\cal U}(\frak l)$.
Thus we can decompose the vector space 
$${\cal U} (\frak l) e_j/R_j = {\cal U} (\frak n_J^+) e_j/ 
\coprod_{i \in J} {\cal U} (\frak n_J^+) (\ad e_i)^{1- 2a_{ij}/a_{ii}}e_j
\oplus {\cal K}e_j.$$

Applying the elimination theorem once again,
using the above decomposition we obtain:
$$L( \coprod_{j \in I\backslash J}{\cal U} (\frak l) e_j/R_j )
= L\l(\coprod_{j \in I\backslash J}\l({\cal U} (\frak n_J^+) e_j/
\coprod_{i \in J} {\cal U} (\frak n_J^+) (\ad e_i)^{1- 2a_{ij}/a_{ii}}e_j\r)\r)
\ltimes \frak c$$
where $\frak c$ is the ideal in $L( \coprod_{j \in I \backslash J}
{\cal U} (\frak l) e_j/R_j )$
generated by the sum of the ${\cal K} e_j$. Each $\frak h^e$-module
${\cal U} (\frak n_J^+) e_j/
\coprod_{i \in J} {\cal U} (\frak n_J^+) (\ad e_i)^{1- 2a_{ij}/a_{ii}}e_j$
is an integrable lowest weight module for the Lie algebra $\frak g_J$, 
denoted by ${\cal U} (\frak n_J^+) \cdot e_j$, with lowest weight $\alpha_j$.
Thus we have a decomposition into semidirect 
products:
$$L(\{e_{i}\}_{i \in I})= L(\{e_{i}\}_{i \in J}) \ltimes 
\l[ \l(L( \coprod_{i \in I \backslash J} {\cal U} (\frak n_J^+)\cdot e_j) 
\ltimes \frak c \r) \ltimes \frak b \r].$$

It is clear that, as ideals of $L(\{e_{i}\}_{i \in I})$, $ \frak b , 
\frak c \subset \frak k_0^+$, and 
$\frak k_J^+$ is $ \frak k_0^+ \cap L(\{e_{i}\}_{i \in J})$. 
Therefore, since all elements of $\frak k_0^+$ are zero in 
$L( \coprod_{i \in I \backslash J} {\cal U} (\frak n_J^+)\cdot e_j)$,
$$L(\{e_{i}\}_{i \in I})/\frak k_0^+ 
= L(\{e_{i}\}_{i \in J})/\frak k_J^+ \ltimes
 L(\coprod_{j \in I \backslash J}{\cal U} (\frak n_J^+)\cdot e_j ). $$
By the definition of $\frak g_J$, $\frak n^+_J =L(\{e_{i}\}_{i \in
J})/\frak k_J^+ $. 
$\square$

\begin{corollary} Let $A$ be a matrix satisfying conditions {\bf C1-C3}.
Assume that the matrix $A$ has only one positive diagonal entry,
$a_{ii} >0$, and if $a_{mj} =0$ then $m=i$, or $j=i$ or $m =j$.  
Let $S =\{ (\adj e_{i})^l e_{j}\}_{0 \leq l \leq -2a_{ij}/a_{ii}}$. 
The subalgebra $\frak n^+ \subset \frak g(A)$ is the
semidirect product of a one-dimensional Lie algebra and a free Lie
algebra,  
$\frak n^+ = \Bbb R e_{i}\oplus L(S).$ Similarly, $\frak n^- = \Bbb R
f_{i}\oplus L(\eta(S)).$ Thus
$$\frak g(A) =  L(S) \oplus (\frak s \frak l_2 + \frak h) \oplus L(\eta(S)).$$ 
\label{thm:fre}
\end{corollary}

The root grading is a grading of the type 
considered in Proposition 5.1 because we have the correspondence
$\alpha_i \mapsto (0,\ldots,1,\ldots,0)$, where $1$ appears in the
$i^{th}$ place. This is the grading (\ref{eq:grad}). The denominator
formula given in the next result is the same as (\ref{eq:denom}), after
the change of variables $e^{\alpha} \mapsto T^{\alpha}$. 

\begin{corollary}
Let $A$ be as in Theorem~\ref{thm:free}, $n_\alpha$ and $T^{\alpha}$ as
in Proposition 5.1. Let $\Delta^J_+ = \Delta \cap \sum_{i \in J}
\Bbb Z_+ \alpha_i$.  The Lie algebra $\frak g(A)$ has 
denominator formula given by (\ref{eq:wil}) and the denominator
formula for $\frak g_J$: 
\begin{eqnarray*}
\prod_{\varphi \in \Delta_+} &&(1-T^\varphi)^{\dim \frak g^{\varphi}}\\
&=& \prod_{\varphi \in \Delta_+^J} (1-T^\varphi)^{\dim 
\frak g_J^{\varphi}} \cdot
\prod_{\varphi \in \Delta_+ \backslash \Delta_+^J} (1-T^\varphi)^{\dim
L^{\varphi}(\coprod_{j \in I \backslash J}{\cal U} 
(\frak n_J^+)\cdot e_j)}
\end{eqnarray*}
$$= \l(\sum_{w \in W} (-1)^{l(w)} T^{w\rho - \rho}\r)\l(
1-\sum_{\varphi \in \Delta^+ \backslash \Delta_+^J}
n_{\varphi}T^{\varphi}\r).$$  \label{cor:den}
\end{corollary}

Obtaining the denominator formula in this way provides an alternative
proof that the radical of the Lie algebra $\frak g$ associated to the
matrix $A$ is zero. This is because Corollary \ref{cor:den}
 uses only the description
of $\frak g$ in terms of generators and relations, while the previous proof
of the denominator formula (see \cite{Jur} or \cite{Kac}) is valid
after the radical of the Lie algebra has been factored out.
Thus we have shown, for the particular type of matrix in
Theorem~\ref{thm:free}:

\begin{corollary}
Let $A$ be as in Theorem~\ref{thm:free}. The generalized Kac-Moody
algebra $\frak g(A)$ has zero radical.
\end{corollary}

\noindent{\bf Remark}: A proof of the fact that any generalized
Kac-Moody algebra has zero radical can by found in \cite{HMY} or
\cite{Jur}. In both cases the argument of \cite{GK} or \cite{Kac} is
extended to 
include generalized Kac-Moody algebras by making use of a lemma
appearing in \cite{Bor1}.

\noindent{\bf Remark}: If we apply (\ref{eq:dim}) to 
$L(\coprod_{j \in I\backslash J} {\cal U} (\frak n_J^+)\cdot e_j)$ we
obtain Kang's \cite{Kan} multiplicity formulas for the special case
of generalized
Kac-Moody algebras with no mutually orthogonal imaginary simple
roots.

\subsection{A Lie algebra related to the modular function $j$}

We will apply our results to an important
example of a generalized Kac-Moody algebra $\frak g(M)$, defined
in Section 6.2 below. 
It will be shown that, if $\frak c$ denotes the center of
the Lie algebra $\frak g(M)$, then $\frak g(M)/\frak c$ is the Monster
Lie algebra.

Recall that the modular function $j$ has the expansion $j(q)=
\sum_{i \in \Bbb Z} c(i) q^i$, where $c(i) =0$ if $i <-1$, $c(-1) =1$,
$c(0) = 744$, $c(1) = 196884$. Let $J(v) = \sum_{i \geq -1} c(i)v^i$ be
the formal Laurent series associated to $j(q) - 744$.

Let $M$ be the symmetric matrix of blocks indexed by
$\{-1,1,2,\ldots\}$, where the block in position $(i,j)$ has entries
$-(i+j)$ and size $c(i) \times 
c(j)$. Thus

$${M} = \left( \begin{array}{c|c|c|c}
          2 & \begin{array}{ccc}  
          0&\cdots&0    \end{array} 
            & \begin{array}{ccc}   -1 & \cdots & -1  \end{array}& \cdots \\
         \hline 
         \begin{array}{c}  0 \\ \vdots \\ 0  \end{array}
                      & \begin{array}{ccc} -2 & \cdots &-2 \\
                                         \vdots & \ddots & \vdots \\
                                       -2 & \cdots & -2   \end{array} &
                      \begin{array}{ccc} -3 & \cdots &-3 \\
                                         \vdots & \ddots & \vdots \\
                                       -3 & \cdots & -3 
                                        \end{array} & \cdots \\
          \hline
          \begin{array}{c}  -1 \\ \vdots \\ -1  \end{array} & 
          \begin{array}{ccc} -3 & \cdots &-3 \\
                                         \vdots & \ddots & \vdots \\
                                       -3 & \cdots & -3 
                                        \end{array} &
          \begin{array}{ccc} -4 & \cdots &-4 \\
                                         \vdots & \ddots & \vdots \\
                                       -4 & \cdots & -4 
                                        \end{array} & \cdots \\
          \hline
          \vdots & \vdots & \vdots 
         \end{array}   \right)     $$

\noindent{\bf Definition 3}: Let $\frak g(M)$ be the generalized
Kac-Moody algebra associated to the matrix $M$ given above. 

We have the standard decomposition
\begin{equation}
\frak g(M)=\frak n^+ \oplus \frak h \oplus \frak n^- .\label{eq:ns}
\end{equation}

The generators of $\frak g(M)$ will be written $e_{jk},f_{jk},h_{jk}$,
indexed by integers $j,k$ 
where $j \in \{-1\} \cup \Bbb Z_+$ and $ 1\leq k \leq c(j)$. Since
there is only one $ e,f$ or $h$ with $j = -1$, we will write these
elements as $e_{-1}, f_{-1}, h_{-1}$. From the
construction of $\frak g(M)$ we see
that the simple roots in $(\frak h^e)^*$ are $\alpha_{-1},
\alpha_{11} ,\alpha_{12}, \cdots , \alpha_{1c(1)},\alpha_{21},
\cdots, \alpha_{2{c(2)}},$ etc. 
Note that for fixed $i$ the functionals
$\alpha_{ij}$ and 
$\alpha_{ik}$ agree on all of $\frak h$ for $1 \leq j,k \leq c(i)$. 
The $\alpha_{ik}$ for
$i > -1$ are simple imaginary 
roots, and the root $\alpha_{-1}$ is the one real simple root. 

\noindent{\bf Remark}: We explain the relationship between our
definition of simple root and that appearing in \cite{Bor3}. If
the restrictions of the simple roots to $\frak h$ are denoted
$\alpha_{-1},\alpha_{1},\alpha_{2},\ldots$, these elements of 
$\frak h^*$ correspond to the notion of ``simple imaginary roots of
multiplicity greater than one'' in \cite{Bor3}. The ``simple root''
$\alpha_i$ has ``multiplicity'' $c(i)$ in Borcherds' terminology.
Fortunately, in this
case, nonsimple roots $\alpha \in (\frak h^e)^*$ do not restrict to
any $\alpha_i$, also ``multiplicities'' do not become infinite,
``roots'' remain either positive or negative, etc. The
functionals $\alpha_i$ are linearly dependent, in fact they span a
two-dimensional space. The
root lattice is described in \cite{Bor3} as the lattice $\Bbb Z
\oplus \Bbb Z$ with the inner product given by 
$\l(\begin{array}{cc} 
      0 & -1 \\
      -1 & 0 \end{array}\r)$.
The ``simple roots'' of this Lie algebra are denoted $(1,n)$ where
$n=-1$ or $n \in \Bbb Z_+$. However, in most cases (including the case
of $\widehat {\frak s \frak l_2}$) serious problems arise when we do not
work in a sufficiently large Cartan subalgebra. (We do not wish to
write down a ``denominator identity'' where some of the terms are
$\infty$). Here, we always work
in $(\frak h^e)^*$, taking specializations when they are illuminating,
as in the case of the denominator identity for $\frak g(M)$ given
below.

Corollary~\ref{thm:fre} applied to the Lie algebra $\frak g(M)$ gives the
following:
\begin{theorem} The subalgebra $\frak n^+ \subset \frak g(M)$ is the
semidirect product of a one-dimensional Lie algebra and a free Lie algebra,
so $\frak n^+ = \Bbb R e_{-1}\oplus L(S).$ Similarly, $\frak n^- = 
\Bbb R f_{-1} \oplus L(S')$. Hence
$$\frak g(M) = L(S) \oplus (\frak s \frak l_2 + \frak h) \oplus L(S').$$
Here $S =\cup_{j \in \Bbb N}
\{ (\adj e_{-1})^l e_{jk}\ |\ 0 \leq l<j, 1 \leq k \leq
c(j)\}$, and $S' = \eta(S)$.
\label{thm:fre2} \end{theorem} $\square$


Now that we have established that $\frak n^+$ is the direct sum of a
one dimensional space and an ideal isomorphic to a 
free Lie algebra we shall obtain the
denominator formula for the Lie algebra $\frak g(M)$.

\begin{corollary}
The denominator formula for the Lie algebra $\frak g(M)$ is
\begin{eqnarray}
\prod_{\varphi \in \Delta_+} (1-T^\varphi)^{\dim \frak g^{\varphi}} 
&=& (1-T^{\alpha_{-1}})\prod_{\varphi \in \Delta_+ \backslash \{\alpha_{-1}\}}
(1-T^\varphi)^{\dim L^{\varphi}(S)}\nonumber \\ 
&=& (1-T^{\alpha_{-1}})\big(1- \sum_{{j \in \Bbb Z_+ \atop 1 \leq k \leq
c(j),\, 0 \leq l <j}} T^{l \alpha_{-1} +\alpha_{jk}}\big),
\end{eqnarray}
which has specialization
\begin{equation}
u(J(u) -J(v)) = \prod_{ i\in \Bbb Z_+ \atop j\in {\Bbb Z_+ \cup \{-1\}} }
  (1-u^iv^j)^{c(ij)} \label{eq:prod} \end{equation}
under the map $\phi :\Delta \rightarrow \Bbb Z \times \Bbb Z$
determined by $\alpha_{ik} \mapsto (1,i)$, where we write $T^{\alpha_{ik}}
\mapsto uv^i$. 
\label{cor1}\end{corollary}

\noindent{\em Proof of Corollary~\ref{cor1}}:

For the denominator identity simply apply Corollary \ref{cor:den}
to the root grading. The $\Bbb Z \times \Bbb Z$-grading of $L(S)$
given above is such that a
generator $(\ad e_{-1})^l e_{j_k}$ has degree $l(1,-1)+(1,j)
=(l+1,j-l)$ with $l <j$. The number of generators of degree $(i,j)$ is
$c(i+j-1)$. 
Applying equation (\ref{eq:wil}) (which is the same as specializing
the second product of the denominator identity via $T^{\alpha_i} \mapsto
uv^i$) gives the formula: 
$$ 1- \sum_{(i,j) \in \Bbb N^2 \backslash\{0\} } c(i+j-1) u^i v^j
  = \prod_{(i,j) \in \Bbb N^2\backslash \{0\}} (1- u^i v^j)^{\dim
L^{(i,j)}}.$$
To obtain the specialization of the denominator formula of $\frak g$
we must include the 
degree $(1,-1)$ subspace $\frak g^{(1,-1)}= \Bbb R e_{-1}$, which is one-%
dimensional. 
$$\prod_{(i,j)} (1-u^i v^j )^{\dim \frak g^{(i,j)}}
=  \prod_{(i,j)\in \Bbb N^2 -\{0\}} (1-u^i v^j )^{\dim L^{(i,j)}} (1- u/v)$$
 $$= (1- \sum_{(i,j) \in \Bbb N^2-\{0\}} c(i+j-1) u^i
    v^j)(1-u/v)$$
$$= 1 - \sum c(i+j-1)u^iv^j - u/v +\sum c(i+j-1)u^{i+1}v^{j-1}$$  
$$  = u(J(u)-J(v)).$$
There is a product formula for the modular function $j$ (see \cite{Bor3})
which can be written:
\begin{equation}
p(j(p)-j(q)) = \prod_{ i= 1,2, \ldots \atop j = -1,1,\ldots}
(1-p^iq^j)^{c(ij)}, \end{equation}
which converges on an open set in $\Bbb C$, and so
implies the corresponding identity for formal power series.
Now we conclude that $\dim \frak g^{(i,j)}
= c(ij)$. $\square$

Note that here we must know the number
theory identity of \cite{Bor3}, to determine the dimension of the
root spaces of $\frak g$. 

The identity (\ref{eq:prod}) is the specialization $e^{\alpha_i}
\mapsto uv^i$ of the denominator identity as it appears when we apply
equation (\ref{eq:denom}) to $\frak g(M)$.

\noindent{\bf Remark}: The matrix $M$ can be replaced by any symmetric
matrix with the same first row (and column) as $M$ with all remaining
blocks having entries strictly less than zero, as long as the minor
obtained by removing the first row and column has the same rank as the
corresponding minor of $M$. 

\noindent{\bf Remark}: Now we apply equation (\ref{eq:dim}) to the
$\Bbb N \times \Bbb N$-grading, and the Lie algebra $L(S)$, where $S
=\{ (\ad e_{-1})^l e_{jk} \ | \ 0 \leq l<j, 1 \leq k \leq c(j)\}$ and
there are $c(i+j-1)$ generators of degree $(i,j)$. Since we already
know the dimension of $L^{(i,j)}$ is $c(ij)$, 
we recover (see \cite{Kan}) the following relations between the
coefficients of $j$: 
\begin{equation}
c(ij) = \sum_{k \in {\Bbb Z}_+ \atop k(m,n) = (i,j)} {1\over k} \mu(k) \sum_{a \in
P(m,n)} {(\sum a_{rs} -1)! \over \prod a_{rs}!}
\prod c(r+s-1)^{a_{rs}}.
\end{equation}

\section{The Monster Lie algebra}
\subsection{Vertex operator algebras and vertex algebras}

For a detailed discussion of vertex operator algebras and vertex
algebras the reader should consult \cite{DL}, \cite{FHL}, \cite{FLM2}
 and the announcement \cite{Bor0}. Results stated here
without proof can either be found in \cite{DL}, \cite{FHL} and
\cite{FLM2} or follow without too much difficulty from the results
appearing there. 

\noindent{\bf Definition 4}: A {\it vertex operator algebra}, $(V,Y,{\bf
1},\omega)$, consists of a vector space $V$, distinguished vectors called 
the {\it vacuum vector} $\bf 1$ and the {\em conformal vector}
$\omega$, and a linear map $Y(\cdot, z): V 
\rightarrow (\mbox{End }V)[[z,z^{-1}]]$ which is a generating function
for operators $v_n$, i.e., for $v \in V,\ Y(v,z)= \sum_{n \in \Bbb Z}
v_n z^{-n-1}$, satisfying the following conditions:
\begin{description}
\item[(V1)] $V =\coprod_{n \in \Bbb Z} V_{(n)}$; for $v \in V_{(n)}$, $n= \wt
(v)$
\item[(V2)] $\dim V_{(n)} < \infty$ for $n \in \Bbb Z$
\item[(V3)] $V_{(n)}=0 $ for $n$ sufficiently small
\item[(V4)] If $u,v \in V$ then $u_nv =0$ for $n$ sufficiently large
\item[(V5)] $Y({\bf 1},z) = 1$
\item[(V6)] $Y(v,z){\bf 1} \in V[[z]]$ and $\lim_{z \rightarrow 0}
Y(v,z){\bf 1}= v$, i.e., the {\it creation property} holds
\item[(V7)] The following {\it Jacobi identity} holds:
$$ z_0^{-1} \delta \l( {z_1-z_2 \over z_0}\r)Y(u, z_1)Y(v,z_2) 
 -z_0^{-1}\delta \l( {z_2-z_1\over -z_0}\r)Y(v,z_2)Y(u,z_1).$$ 
\begin{equation}
 = z_2^{-1} \delta \l({z_1-z_0 \over z_2}\r) Y(Y(u,z_0)v,z_2)
\end{equation}
\end{description}
The following conditions relating to the vector $\omega$ also hold;
\begin{description}
\item[(V8)] The operators $\omega_n$ generate a Virasoro algebra i.e., if
we let $L(n) = \omega_{n+1}$ for $n \in \Bbb Z$ then
\begin{equation}
[L(m),L(n)] = (m-n)L(m+n) + (1/12) (m^3-m)\delta_{m+n,0}(\mbox{rank} V)
\end{equation}
\item[(V9)] If $v \in V_{(n)}$ then $L(0)v = (\wt v)v = nv$
\item[(V10)] ${d \over dz} Y(v,z)= Y(L(-1)v,z)$.
\end{description}

\noindent{\bf Definition 5}: A {\it vertex algebra} $(V,Y,{\bf 1}, \omega)$
is a vector space $V$ with all of the above properties except for 
$\bf V2$ and $\bf V3$.

\noindent{\bf Remark}: This definition is a variant, with $\omega$, of
Borcherds' original definition of vertex algebra in \cite{Bor0}.

An important class of examples of vertex algebras (and vertex operator
algebras) are those associated with lattices. For the sake of the
reader who may be unfamiliar with the notation we will briefly review 
this construction in the case of an even lattice. For complete details
(and more generality) the reader may consult \cite{FLM2} or \cite{DL}.
Given an even lattice $L$ one can construct a vertex algebra $V_L$ with
underlying vector space:
$$V_L=S(\hat \frak h^-_\Bbb Z)\otimes \Bbb R\{L\}.$$
Here we take $\frak h= L \otimes_{\Bbb Z} \Bbb R$, and 
$\hat \frak h^-_{\Bbb Z}$ is the
negative part of the Heisenberg algebra (with $c$ central) defined by:
$$\hat {\frak h}_{\Bbb Z} =\coprod_{n \in \Bbb Z} \frak h \otimes t^n  
\oplus {\Bbb R} c \subset \frak h \otimes {\Bbb R} [[t]] \oplus 
{\Bbb R} c.$$ 
Therefore, 
$$\hat \frak h^-_\Bbb Z = 
\coprod_{n <0} \frak h \otimes t^n.$$
 The symmetric algebra on $\hat \frak h_{\Bbb Z}^-$ is denoted $S(\hat
\frak h^-_\Bbb Z)$. Given a central extension of $L$ by a group of order $2$
i.e.,
$$ 1 \rightarrow \la \kappa | \kappa^2 =1\ra \rightarrow \hat {L}
 {\buildrel - \over \rightarrow} L \rightarrow 1,$$
with commutator map given by $\kappa^{\la \alpha, \beta\ra}$, $
\alpha, \beta \in L$, $\Bbb R$ is given the structure of a nontrivial
$\la \kappa \ra$-module.
 Define $\Bbb R \{L\}$ to be the induced
representation
$\mbox{Ind}_{\la \kappa \ra}^{\hat L} {\Bbb R}$. 

If $a \in {\hat L}$
denote by $\iota(a)$ the element $a \otimes 1 \in \Bbb R \{L\}$. We
will use the  
notation $\alpha(n) = \alpha \otimes t^n \in S(\hat \frak h^-_\Bbb Z)$.
The vector space $V_L$ is spanned by elements of the form:  
$$\alpha_1(-n_1)\alpha_2(-n_2)\ldots\alpha_k(-n_k)\iota(a)$$
where $n_i \in \Bbb N$.
The space $V_L$, equipped with $Y(v,z)$ as defined in \cite{FLM2}
satisfies properties $\bf V1$ and $\bf V4 - V10$, so
is a vertex algebra with conformal vector $\omega$. Features of
$Y(v,z)$ to keep in mind from \cite{FLM2} are: $\alpha(-1)_n =\alpha(n)
\mbox{ for all } n \in 
\Bbb Z$; the $\alpha(n)$ for $n <0$ act by left multiplication on
$u \in V_L$; and
\[
\alpha(n)\iota(a) = \left\{ \begin{array}{cc}
                         0 & \mbox{ if $n>0$ }\\
                 \la\alpha, \bar a\ra\iota(a) &\mbox{ if $n=0$}.     
                       \end{array} \right. 
\]

\noindent{\bf Definition 7:} A bilinear form $(\cdot, \cdot)$ on a
vertex algebra $V$ is {\em invariant} (in the sense of
\cite{FHL}) if it satisfies 
\begin{equation}
( Y(v,z)w_1,w_2 ) =
( w_1, Y(e^{zL(1)}(-z^{-2})^{L(0)}v,z^{-1})w_2 ) . \label{eq:inv}
\end{equation}
(By definition $x^{L(0)}$ acts
on a homogeneous element $v \in V$ as multiplication by $x^{wt(v)}$).
Such a form satisfies $( u, v )=0$ unless $\wt (u) = \wt (v)$. 

\begin{lemma} Let $L$ be an even unimodular lattice. There is a
nondegenerate symmetric invariant bilinear form $(\cdot,\cdot )$ on
$V_L$.\label{lem:inv1}
\end{lemma}

\noindent{\em Proof}:
The vertex algebra $V_L$ is a module for itself under the adjoint 
action. In fact, $V_L$ is an irreducible module and any irreducible
module of $V_L$ is isomorphic to $V_L$ \cite{D}.

In order to define the contragredient module note that $V_L$
is graded by the lattice $L$ as well as by weights, and that under
this double grading $\dim {V_L}^r_{(n)} < \infty$ for $r \in L, n \in
\Bbb Z$. Let ${V_L}' = \coprod_{n \in \Bbb Z\atop  r \in L}
({V_L}^r_{(n)})^*$, the restricted dual of $V_L$. Denote by
$\langle \cdot, \cdot \rangle$ the natural pairing between
$V_L$ and ${V_L}'$. Results of \cite{FHL} pertaining to adjoint vertex
operators and the contragredient module now apply to
${V_L}'$. In particular, the space ${V_L}'$ can be given the structure of a
$V_L$-module $({V_L}', Y')$ via
$$\langle Y'(v,z)w',w  \rangle
=\langle w',Y(e^{zL(1)}(-z^{-2})^{L(0)}v,z^{-1})w \rangle
$$
for $v,w \in V_L$ $w'\in {V_L}'$.
Since the adjoint module $V_L$ is irreducible, the contragredient
module ${V_L}'$ is also irreducible. By the result of \cite{D} quoted
above, $V_L$ is isomorphic to ${V_L}'$ as a $V_L$-module, which is
equivalent to $V_L$ having a nondegenerate invariant bilinear form.
(See remark 5.3.3 of \cite{FHL})
$\square$

The ``moonshine module'' 
$V^\natural$ is an infinite-dimensional representation of the Monster
simple group constructed and shown to be a vertex
operator algebra in \cite{FLM2}. The graded dimension of $V^\natural$
is $J(q)$. There is a positive definite bilinear form $(\cdot, \cdot)$
on $V^\natural$.
The vertex operator algebra $V^\natural$ satisfies
all of the conditions of the no-ghost theorem (Theorem~\ref{thm:ng}),
with $G$ taken to be the Monster simple group. This
vertex operator algebra will be essential to the construction of the
Monster Lie algebra.

\begin{lemma} The positive definite form $(\cdot, \cdot)$ on
$V^\natural$ defined in \cite{FLM2} is invariant. 
\label{lem:inv2}
\end{lemma}

\noindent{\em Proof}: There is a unique up to constant multiple
nondegenerate symmetric invariant bilinear form on $V^\natural$
\cite{Li2}. Fix such a form $(\cdot,\cdot)_{1}$ by taking
$(1,1)_{1}=1$.

Let $u \in V^\natural_{(2)}$, a homogeneous element of
weight $2$. By invariance $$(u_nw_1,w_2)_{1}=(w_1,u_{-n + 2}w_2)_{1}$$
for $w_1,w_2\in V^\natural$.
We claim
\begin{equation}
(u_nw_1,w_2)=(w_1,u_{-n + 2}w_2)
\label{eqn:grs}
\end{equation}
$w_1,w_2\in V^\natural$. 
In order to prove equation (\ref{eqn:grs}) we recall the construction and
properties of $V^\natural$ of \cite{FLM2}.
Let $x_a^+ = \iota(a) +\iota(a^{-1})$ for $a \in \wedge$ (the Leech
lattice),
$$\frak k =
S^2(\frak h \otimes t^{-1})\oplus
\sum_{a \in \hat{\wedge}_4} \Bbb R x_a^+,$$
and let $\frak p $ be the space of elements of $V^T_\lambda$ (the
``twisted space'') of weight 2. 
 Then $V^\natural_{(2)}= \frak k \oplus \frak p$. The action $Y(v,z)$
of elements $v \in \frak p$ is determined by conjugating by certain
elements of $\Bbb M$ the Monster simple group (see 12.3.8 and 12.3.9 of
\cite{FLM2}). Conjugation by these elements map $v \in\frak p$ to
$\frak k$. Since the form  $(\cdot,\cdot)$ is invariant under $\Bbb
M$, it is sufficient to check that equation (\ref{eqn:grs}) holds for
elements of $\frak k$. Therefore, 
it suffices to check
 equation (\ref{eqn:grs}) for two types of elements
$x_a^+, a \in\hat{\wedge}_4$ 
and $g(-1)^2, g(-1) \in \frak h \otimes t^{-1}$.
For $u = x_a^+$
equation (\ref{eqn:grs}) follows immediately from \cite{FLM2}. For $u=g(-1)^2$
\begin{eqnarray*}
Y(u,z)&=& \mbox{$\circ\atop\circ$} g(z)^2 \mbox{$\circ\atop\circ$}\\
&=& g(z)^- g(z) +g(z)g(z)^+ 
\end{eqnarray*}
Using $(g(i)w_1,w_2) = (w_1,g(-i)w_2)$
one computes the adjoint of $g(z)^- g(z) +g(z)g(z)^+ $ which is
$$g(z^{-1}) g(z^{-1})^+ + g(z^{-1})^- g(z^{-1}) = Y(u,z^{-1})z^{-4}. $$
Thus equation (\ref{eqn:grs}) holds for all $u \in \frak k$, and so
for all $u \in V^\natural_{(2)}$.

Now recall that $V^\natural_{(2)}=\cal{B}$, the Griess algebra, and
the notation $\hat{\cal{B}}$ for the commutative affinization of the
algebra $\cal{B}$,
$$\hat{\cal B} = {\cal B} \otimes {\Bbb R}[t,t^{-1}] \oplus {\Bbb R} e$$
where $t$ is an indeterminate and $e \neq 0$ (with nonassociative
product given in \cite{FLM2}). By Theorem 12.3.1 
\cite{FLM2} $V^\natural$ is an irreducible graded
$\hat{\cal{B}}$-module, under
\begin{eqnarray*}
\pi: \hat{\cal{B}} &\rightarrow & \mbox{End}V^\natural\\
      v\otimes t^n &\mapsto & x_v(n) \ v \in {\cal B}\\
        e &\mapsto& 1
\end{eqnarray*}

Schur's lemma then implies that any nondegenerate symmetric bilinear
form satisfying
equation (\ref{eqn:grs}) is unique up to multiplication by a constant.
Thus we can conclude that $(\cdot,\cdot)_{1} = (\cdot,\cdot)$, since
the length of the vacuum is one with respect to each form.
$\square$

Given $V$ a vertex
operator algebra, or a vertex algebra with $\omega$ 
and therefore an action of the Virasoro algebra, let $$P_{(i)} = \{ v \in V
| L(0)v = iv , L(n)v = 0 \mbox{ if } n>0\}.$$ Thus $P_{(i)}$ consists of
the lowest weight vectors for the Virasoro algebra of
weight $i$. Then $P_{(1)}/L(-1)P_{(0)}$ is a Lie algebra with bracket
given by $[u +L(-1)P_{(0)} ,v+L(-1)P_{(0)}] = u_0v + L(-1)P_{(0)}$.
If the vertex algebra $V$ has an invariant 
bilinear
form $( \cdot, \cdot )$ this induces a form 
$( \cdot, \cdot)_{Lie}$ on the Lie algebra
$P_1/L(-1)P_{(0)}$, because $L_{-1}P_{(0)} \subset
\mbox{rad}(\cdot,\cdot)$. 
Invariance of the form on the vertex algebra implies for 
$u,v \in P_{(1)}$:
\begin{equation}
( u_0 v, w) = - ( v,u_0 w).
\end{equation}
So that the induced form is invariant on the Lie algebra 
$P_{(1)}/L(-1)P_{(0)}$.

Tensor products of vertex operator algebras are again vertex operator
algebras (see \cite{FHL}), and more generally, by \cite{DL}, the tensor
product of vertex algebras is also a vertex
algebra. Given two vertex algebras $(V,Y,{\bf 1}_V, \omega_V)$ and
$(W,Y,{\bf 1}_W, \omega_W)$ the vacuum of
$V \otimes W$ is ${\bf 1}_V \otimes {\bf 1}_W$ 
and the conformal vector $\omega$ is given by $\omega_V \otimes 
{\bf 1}_W + {\bf 1}_V \otimes \omega_W$. 
If the vertex algebras $V$
and $W$ both have invariant forms then it is not difficult to show
that the form on $V \otimes W$ given by the product of the forms on
$V$ and $W$ is also invariant in the sense of equation (\ref{eq:inv}).

\subsection {The Monster Lie algebra }

We will review the construction of the Monster Lie algebra given in
\cite{Bor3}. Then we give a theorem regarding its structure as a quotient
of $\frak g(M)$. Let ${\cal L}=\Bbb Z \oplus \Bbb Z$ with bilinear form
$\la \cdot,\cdot \ra$ given by the matrix 
$\l(\begin{array}{cc}
      0 & -1 \\
      -1 & 0 
\end{array}\r)$. 

\noindent{\bf Remark}: $\cal L$ is the rank two Lorentzian lattice,
denoted in \cite{Bor3} as $II_{1,1}$. 

Fix a symmetric invariant bilinear form $(\cdot,\cdot)$
 on $V_{\cal L}$, normalized by taking $( {\bf 1},{\bf 1} ) =-1$. The
reason that we choose this normalization is so that the resulting
invariant bilinear form on the monster Lie algebra will have the usual
values with respect to the Chevalley generators, and so that the
contravariant bilinear form defined below will be positive definite
and not negative definite on nonzero weight spaces.
Denote by  $(\cdot,\cdot )$ the symmetric invariant 
bilinear form on $V^\natural\otimes V_{\cal L}$ given by the product of
the invariant bilinear forms on $V^\natural$ and $V_{\cal L}$.

\noindent{\bf Definition 7}: The {\it Monster Lie algebra} $\frak m$
is defined by 
$$\frak m
=P_1/\mbox{rad}(\cdot,\cdot)_{Lie} =
(P_1/L_{-1}P_0)/\mbox{rad}(\cdot,\cdot)_{Lie}.$$

When no confusion will arise, we will use the same notation
for the invariant form on the vertex algebra, and for the induced form
on the Lie algebra.

Note that, by invariance
$( e^r,e^s ) = (-1)^{\la r,r\ra /2}\la 1,(e^r)_{(-1 + \la
r,r\ra)}e^s\ra= 0$ unless $r=-s \in \cal L$. Therefore, the induced
form on $\frak m$ satisfies the condition that $\frak m_r$ be
orthogonal to $\frak m_s$ if $r\neq -s \in \cal L$. 

\noindent{\bf Definition 8}: Let $\theta$ be 
the involution of $V_{\cal L}$ given by $\theta \iota (a) =
(-1)^{{ wt} (a)} \iota(a^{-1})$ and $\theta (\alpha(n)) = -\alpha(n)$.
This
induces an involution $\theta$ on all of $V^{\natural} \otimes
V_{\cal L}$ by letting $\theta (u \otimes v e^r) = 
u \otimes \theta(v e^r)$. Use the same notation for the
involution induced by $\theta$ on $\frak m$.

Note that $\theta : \frak m_r
\rightarrow 
\frak m_{-r}$ if $r \neq 0$.

Let $(\cdot, \cdot)_0$ be the contravariant bilinear form
on $V^{\natural} \otimes V_{\cal L}$ given by 
$(u, v)_0 = - (u,\theta (v))$, 
$u,v \in V^{\natural} \otimes V_{\cal L}$.
We also denote by $(\cdot, \cdot)_0$ 
the contravariant bilinear form
on $\frak m$ given by $(u, v)_0 = - (u,\theta (v))_{Lie}$, 
$u,v \in \frak m$.

Elements of $\frak m$ can be written as $\sum u 
\otimes v e^r$, where $u \in V^{\natural}$ and $ ve^r = v \iota (e^r)
 \in V_{\cal L}$. Here, a section of the map $\hat {\cal L}
 {\buildrel - \over \rightarrow} {\cal L}$ has been chosen so that $e^r \in
{\hat {\cal L}}$ satisfies $\overline{e^r}= r 
\in {\cal L} $. There is a grading of $\frak m$ by the lattice defined by
$\deg (u \otimes ve^r)=r$. 

Recall the definition of the Lie algebra $\frak g(M)$ and the standard
decomposition.

\begin{theorem} Let $\frak c$ denote the center of the Lie algebra 
$\frak g(M)$. Then $$\frak g(M)/\frak c = \frak m.$$
There is a  triangular decomposition
$\frak m = \frak m^+ \oplus \frak h \oplus \frak m^-$, where $\frak h
\cong \Bbb R \oplus \Bbb R$. The subalgebras $\frak m^{\pm}$ are
isomorphic to $\frak n^{\pm} \subset \frak g(M)$.
\label{thm:mons} 
\end{theorem}

Theorem~\ref{thm:mons} is proven after the statement of the no-ghost
theorem. In \cite{Bor3} the no-ghost theorem from string
theory is used to see that the Monster Lie algebra has homogeneous
subspaces isomorphic to $V^\natural_{(1+mn)}$. A precise statement
of the no-ghost theorem as it is used here is provided for the reader,
and a proof is given in the appendix.

\begin{theorem}[no-ghost theorem] Let $V$ be a vertex operator algebra
with the following properties:
\begin{enumerate}
 \item [i.] $V$ has a symmetric invariant nondegenerate bilinear form.
 \item [ii.] The central element of the Virasoro algebra acts as
multiplication by 24.
 \item [iii.] The weight grading of $V$ is an $\Bbb N$-grading of $V$,
i.e., $V = \coprod_{n=0}^\infty V_{(n)}$, and $\dim V_{(0)}=1$.
 \item [iv.] $V$ is acted on by a group $G$ preserving the above
structure; in particular the form on $V$ is $G$-invariant. 
\end{enumerate}
Let ${\cal P}_{(1)} = \{ u \in V \otimes V_{\cal L} | L_0 u
=u, L_iu =0, i>0\}$. The group $G$ acts on $V \otimes V_{\cal L}$ via the
trivial action on $V_{\cal L}$. Let ${\cal P}_{(1)}^r$ denote the
subspace of ${\cal P}_{(1)}$ of degree $r \in 
\cal L$. Then the quotient of ${\cal P}_{(1)}^r$ by the nullspace of its
bilinear form is isomorphic as a $G$-module with $G$-invariant
bilinear form to $V_{(1- \la r,r\ra /2)}$ if $r \neq 0$ and to $V_{(1)}
\oplus {\Bbb R^2}$ if $r =0$. \label{thm:ng} \end{theorem}


\noindent{\it Proof of Theorem~\ref{thm:mons}}\/:

The no-ghost theorem applied to $V^\natural$ immediately gives
$\frak m_{(m,n)} \cong V^\natural_{(mn+1)}$ if $(m,n) \neq (0,0)$. Thus
the elements 
of $\frak m_r$ where $r \neq 0$ are spanned by elements of the form $u
\otimes e^r$ where $r \in \cal L$, and $u \in V^\natural$ is
an element of the appropriate weight.
(We will use elements of
$V^\natural \otimes V_L$ to denote their equivalence classes in 
$\frak m$.) 

We will show that all of the conditions of
Theorem~\ref{theorem:hom}  
are satisfied. 
By considering the weights (with respect to $L_0$), we see that
the abelian subalgebra $\frak m_{(0,0)}$ is spanned by elements of the
form $ 1 \otimes  \alpha(-1)\iota(1)$ where $\alpha \in {\cal L}
\otimes_{\Bbb Z}\Bbb R$. Note that $\frak m_{(0,0)}$ is 
two-dimensional. 

By definition $\theta= -1$ on $\frak m_{(0,0)}$. 
Thus $(\theta (x),\theta (y))= ( x,y)$ for
$x,y \in \frak m_{(0,0)}$. It 
also follows from the definition, and symmetry of the form on $V_L$
that 
\begin{eqnarray*}
(\theta (u \otimes e^r),\theta (v \otimes e^{-r}))
&=& ( u\otimes e^{-r},v\otimes e^r) \\
&=& ( u\otimes  e^r,v\otimes e^{-r})
\end{eqnarray*}
for $u, v \in V^\natural, r \in
{\cal L}$ .

Consider 
$-( x, \theta (x) )$ for $x \in\frak m_r, r \neq 0$. To 
see that this is strictly positive it is enough to consider
elements of the form
$x = u \otimes e^r$ where $u \in V^\natural$. Recalling the
normalization of the form on $V_L$,
\begin{eqnarray*}
\lefteqn{ ( x,\theta (x) ) } \\  
&&= ( u,u ) (-1)^{\la r,r\ra/2} ( e^r ,e^{-r} ) \\
&&= (-1)^{\la r,r \ra } ( u,u ) 
( {\bf 1}, (e^r)_{(-1 + \la r,r\ra)}e^{-r})\\
&& =  ( u,u )({\bf 1},{\bf 1}) <0.
\end{eqnarray*}
Therefore, we have the desired 
properties on the form $(\cdot,\cdot)=(\cdot,\cdot)_{Lie}$, and the
contravariant form $(\cdot,\cdot)_0$.

Now if we grade 
$\frak m$, as in \cite{Bor3}, by $ i= 2m+n \in\Bbb Z$ then we see that $\frak
m$ satisfies the grading condition. Furthermore, condition {\em 3}
is satisfied if we take $ \theta$ to be the involution.

Let $v \in V^{\natural}$, so
that $v \otimes  e^{r}$ is of degree $ r=(m,n)$. Then
\begin{equation}
 (1 \otimes\alpha(-1)\iota(1))_0 v\otimes e^{r} =v \otimes
\alpha(0)  e^{r}= \langle\alpha,r\rangle v\otimes e^{r} .
\end{equation}
Thus $1 \otimes \alpha(-1)\iota(1)$ acts as the scalar
$\langle\alpha,r\rangle$ on $\frak m_{(m,n)}$. 
Thus all elements of $\frak m_{(0,0)}$
act as scalars 
on the $\frak m_{(m,n)}$. As $\alpha$ ranges over
$\Bbb R\oplus \Bbb R$ the action distinguishes between spaces of
different degree. 
This  establishes conditions {\em 1, 2} and {\em 3} of
Theorem~\ref{theorem:hom}.

To see that 
$\frak m_{(0,0)} \subset 
[\frak m, \frak m ]$, let $u,v \in V_{(2)}^\natural$ and $a = e^{(1,1)}$, 
$b =e^{(1,-1)}$. We will show
\begin{equation}
[u \otimes\iota(a), v \otimes \iota(a^{-1})] \label{eq:g}
\end{equation}
 and 
\begin{equation}
[\iota(b),\iota(b^{-1})] \label{eq:be}
\end{equation}
 are two linearly independent vectors in
$\frak m_{(0,0)}$. Since we know that $\frak m_{(0,0)}$ is two-%
dimensional, this will give condition {\em 4} of
Theorem~\ref{theorem:hom}.

By \cite[8.5.44]{FLM2} we have that formula (\ref{eq:be}) is $\iota(b)_0
\iota(b^{-1})= \bar b(-1)\iota(1)$.
The formula \cite[8.5.44]{FLM2} also shows $\iota(a)_i \iota(a_{-1}) =
0$ unless $i \leq -3$ and
$$\iota(a)_i \iota(a^{-1}) =\l\{\begin{array}{cc}
                                \iota(1) & \mbox{ if $i =-3$}\\
                                {\bar a}(-1)\iota(1) & \mbox{ if $i =-4$}.
                                \end{array} \r. $$
Using the Jacobi identity (19) or its component form
\cite[8.8.41]{FLM2} and the above, we obtain for formula (\ref{eq:g}):
\begin{eqnarray*}
(u \otimes\iota(a))_0 (v \otimes \iota(a^{-1})) \\
=  u_3v \otimes{\bar a}(-1)\iota(1)&& \\
= c {\bf 1}\otimes{\bar a}(-1)\iota(1). &&
\end{eqnarray*}
Since we can pick $u$ and $v$ such that $c \neq 0$ these vectors are
linearly independent and we are done. 

By definition the radical of the bilinear 
form on $\frak m$ is zero, so by Corollary 4.1, $\frak m$ is 
$\frak l/\frak c $ for some generalized  
Kac-Moody algebra $\frak l$. In fact  
 $\frak m = \frak g(M)/ \frak c$:
Because $\frak m_{(0,0)}$ is the image 
of a maximal toral subalgebra, it must also be a maximal toral 
subalgebra. Define {\em roots} of $\frak m$ as elements
$\alpha \in (\frak m_{(0,0)})^*$ such that $ [h, x] = \alpha(h) x$ for
all $x \in \frak m$. 
The grading given by the lattice $\cal L$
corresponds to the 
root grading of $\frak m$ because we have shown that the elements of
$\frak m_{(0,0)}$ act as scalars 
on the $\frak m_{(m,n)}$, and that the spaces of different degree are
distinguished. It follows from the no-ghost
theorem applied to $V^\natural$ that  
$\frak m_{(m,n)} \cong V^\natural_{mn+1}$. We know from \cite{FLM2} that
$\dim V^\natural_{mn+1} = c(mn)$. 
Consider the roots of $\frak g(M)$ restricted to $\frak h$. Then the
dimensions of these restricted root spaces of $\frak g(M)$ are given by 
$c(mn)$. By the specialization (12) of the denominator formula for 
$\frak g(M)$ the generalized Kac-Moody algebra $\frak g(M)/\frak c$
is isomorphic to $\frak m$. $\square$

Since the map  given by Corollary 4.1 is an isomorphism on 
$\frak n^{\pm}$ there are immediate corollaries to
Theorem~\ref{thm:fre2} and the denominator identity for $\frak g(M)$.
\begin{corollary}
The Monster Lie algebra $\frak m$ can be written as $\frak m = 
\frak u^+ \oplus \frak g\frak l_2 \oplus \frak u^-$, where $\frak u^\pm$
are free Lie algebras with countably many generators given by
Corollary 5.1.
\end{corollary}

\begin{corollary}
The Monster Lie algebra has the denominator formula:
\begin{equation}
u(J(u) -J(v)) = 
\prod_{ i \in \Bbb Z_+\atop  j\in {\Bbb Z_+ \cup \{-1\}}}
   (1-u^iv^j)^{c(ij)}.
\end{equation}
\end{corollary}

\vfill\eject
\appendix
\section{The Proof of the no-ghost theorem}

In \cite{Bor3} it is shown how to use the no-ghost theorem from string
theory to understand some of the structure of the monster Lie algebra.
The proof of that theorem, Theorem~\ref{thm:ng}, is reproduced here 
with the necessary
rigor and in a more algebraic context.

 The space $V \otimes V_{\cal L}$ is a vertex
algebra with conformal vector.
Recall from Section 6.1 that elements of the Virasoro algebra 
acting on $V \otimes V_{\cal L}$ 
satisfy the relations: 
$$[L_i, L_j]=(i-j)L_{i+j} + {26 \over 12}(i^3 - i)\delta_{i+j,0}.$$
Given a nonzero $r \in {\cal L}$, fix nonzero $w \in \cal L$ 
such that $\la w,w \ra =0$ and $\la r,w
\ra \neq 0$. Define operators $K_i$, $i \in {\Bbb Z}$, on 
$V \otimes V_{\cal L}$ by $K_i
=({ 1}\otimes w(-1))_{i}= { 1}\otimes w(i)$. 

Let $\cal A$ be the Lie algebra generated by the operators $ L_i, K_i$ with
$i \in \Bbb Z$. These operators satisfy the relations:
$$\displaylines {[L_i, K_j] = -jK_{i+j} \cr [K_i, K_j]=0.}$$
The first relation follows from the formula (in $V_{\cal L}$) $[L_m,
w(n)] = -n(w(n+m))$ of \cite[8.7.13]{FLM2} and the fact that $[L_m
\otimes { 1}, { 1} \otimes w(n)] =0$. The second relation holds
because $ [w(i),w(j)] = \la w,w \ra i \delta_{i+j,0} $  \cite[8.6.42]{FLM2}
and $\la w,w\ra =0$. 

Let ${\cal W}$ be the Virasoro subalgebra of $\cal A$ generated by the
$ L_i, i \in \Bbb Z$, and let $\cal Y$ be the abelian subalgebra generated
by the $K_i, i \in \Bbb Z$. 
Denote by ${\cal A}^+$ the subalgebra generated by the $ L_i, K_i$ with $i>0$,
let ${\cal A}^-$ be the subalgebra generated by the $ L_i, K_i$ 
with $i < 0$, and let ${\cal A}^0$ be the subalgebra generated by
$L_0$, $K_0$.
The subalgebras ${\cal W}^\pm$, ${\cal W}^0$ and $\cal Y^\pm$, 
${\cal Y}^0$ are defined analogously.

The vertex algebra $ V \otimes V_{\cal L}$ is graded by $\cal L$,
because $V_{\cal L}= S(\hat \frak h^-_{\Bbb Z}) \otimes 
{\Bbb R}\{\cal L\} $ has
such a grading. The subspace of degree $r$ is $V \otimes S(
\hat \frak h^-_{\Bbb Z}) \otimes e^r$. This space will be denoted
$\cal H$. The following subspaces of the ${\cal A}$-module
$\cal H$ will be useful:
$$\displaylines { {\cal P}= \{v \in {\cal H}| {\cal W}^+v =0\},  \cr
                  {\cal T} = \{v \in {\cal H}| {\cal A}^+v=0\}, \cr
                  {\cal N} = \hbox {the radical of the bilinear 
                   form $(\cdot, \cdot)_0$
                   on $\cal P$ },\cr
                  {\cal K} = U(\cal Y)  {\cal T}. }$$
Denote $V \otimes e^r $ by  $ Ve^r$.

\begin{lemma} With respect to the bilinear form $( \cdot, \cdot)_0$ on
$V \otimes V_{\cal L}$, $L_i^* = L_{-i}$
and $K_i^* =K_{-i}$ for all $i \in \Bbb Z$.\label{lem:dual}\end{lemma} 

\noindent {\it Proof}\/: Let $\omega$ be the conformal vector of 
$V \otimes V_{\cal L}$, so $L_i = \omega_{i+1}$ and $\theta(\omega)=\omega$.
By the definition of the form $( \cdot,\cdot )_0$ and 
equation (\ref{eq:inv}) $L_i^*$ is 
$$\mbox{Res}_{z^{-i-2}} 
Y(e^{zL_1}(-z^{-2})^{L_0}\omega, z^{-1}).$$ Since $L_1\omega =0$ and
$\wt \omega =2$ we have 
$$Y(e^{zL_1}(-z^{-2})^{L_0}\omega, z^{-1})
= \sum_{n \in \Bbb Z} \omega_nz^{n-3}.$$
Thus $L_i^* = \omega_{-i+1}= L_{-i}$.

Now consider $K_i = (1 \otimes w(-1))_i$. By equation (\ref{eq:inv})
$$K_i^* = \mbox{Res}_{z^{-i-1}} 
\theta Y(e^{zL_1}(-z^{-2})^{L_0}(1 \otimes w(-1)), z^{-1}).$$ To calculate 
this, note that
$\theta(1 \otimes w(-1))= -(1 \otimes w(-1))$, that 
$\wt (1 \otimes w(-1))=1$ and that $L_1 (1 \otimes w(-1))=0$, so 
$$-Y(e^{zL_1}(-z^{-2})^{L_0}(1 \otimes w(-1)), z^{-1})= 
\sum_{n\in \Bbb Z}(1 \otimes w(-1))_n z^{n-1}.
$$
We conclude $K_i^* = K_{-i}$.$\square$

\begin{lemma} The bilinear form $( \cdot, \cdot)_0$ restricted to 
${\cal H}$ is nondegenerate.
\label{lem:nondeg}
\end{lemma} 

\noindent {\it Proof}\/:
The form $( \cdot, \cdot)_0$ on $V \otimes V_{\cal L}$ is nondegenerate.
The form also satisfies $( u, v)_0=0$ unless $\deg(u)=\deg(v)$ in $\cal L$.
Thus the radical of the form $( \cdot, \cdot)_0$ restricted to 
${\cal H}$ is contained in the 
radical of the form on $V \otimes V_{\cal L}$.$\square$

\begin{lemma} ${\cal H} = U({\cal A})  {\cal T}$. 
\label{lem:HAT}
\end{lemma} 

\noindent {\it Proof}\/:
The bilinear form on ${\cal H}$ is nondegenerate (Lemma
\ref{lem:nondeg}) and distinct $L_0$-weight spaces of $\cal H$ are
orthogonal.
Thus the finite-dimensional $i$th $L_0$-weight space 
${\cal H}_i = U({\cal A})  {\cal T}_i \oplus ( U({\cal A})  {\cal T})_i^\perp.$
Then
there is a decomposition into ${\cal A}$-submodules:
$${\cal H}= U({\cal A})  {\cal T} \oplus ( U({\cal A})  {\cal T})^\perp.$$ 
If the graded submodule $( U({\cal A})  {\cal T})^\perp$ is nonempty
then it contains a vector
annihilated by ${\cal A}^+$ by the following argument: The grading of ${\cal H}$
(by weights of $L_0$) is such that ${\cal H}= 
\coprod_{i \geq {1\over 2}\la r,r \ra}{\cal H}_i$.
The actions of the generators $L_i$ and $K_i$, $i >0$, of ${\cal A}^+$
lower the weight of a vector in ${\cal H}$. If
$n$ is the smallest integer such that $( U({\cal A})  {\cal T})^\perp \cap
{\cal H}_n$ is nonzero, then this subspace consists of vectors 
annihilated by ${\cal A}^+$. By definition such a vector is in $\cal T$, hence
is in $U({\cal A})  {\cal T}$, a contradiction.
$\square$

\begin{lemma}
${\cal K} = {\cal T}\oplus \mbox{rad}(\cdot,\cdot)_0$
\label{lem:decomp1}
\end{lemma}

\noindent {\it Proof}\/:
Note that ${\cal Y}= {\cal Y}^- \oplus {\cal Y}^+ \oplus {\cal Y}^0$,
so that by the  
Poincare-Birkhoff-Witt
theorem ${\cal K} = U({\cal Y}^-)U({\cal Y}^+) U({\cal Y}^0){\cal T}$. 
By definition of 
${\cal T}$,
$U({\cal Y}^+) {\cal T}={\cal T}$ so ${\cal K} = U({\cal Y}^-) {\cal T}$. 
Thus 
$${\cal K} = {\cal T}\oplus {\cal Y}^- U({\cal Y}^-) {\cal T}.$$
By Lemma \ref{lem:dual}
\begin{eqnarray*}
( {\cal K}, {\cal Y}^- U({\cal Y}^-) {\cal T})_0 
&=& ( {\cal Y}^+{\cal K},U({\cal Y}^-) {\cal T})_0 \\
&=&0
\end{eqnarray*}
Therefore
 ${\cal Y}^- U({\cal Y}^-) {\cal T}\subset \mbox{rad}(\cdot,\cdot)_0$.
Furthermore, ${\cal T}\cap\mbox{rad}(\cdot,\cdot)_0=0$ because if 
$t \in {\cal T}\cap\mbox{rad}(\cdot,\cdot)_0$ then
$( t, U({{\cal A}}){\cal T})_0=( t, U({\cal A}^-){\cal T})_0=0$, and by Lemma
\ref{lem:nondeg} $t =0$. $\square$ 

\begin{lemma}
${\cal K} = Ve^r \oplus \mbox{rad}(\cdot,\cdot)_0$
\label{lem:decomp2}
\end{lemma}

\noindent {\it Proof}\/:
It is immediate from the definition that elements of ${\cal K}$ 
are lowest weight vectors of ${\cal Y}$. Furthermore,
\begin{eqnarray*}
\cal H &=& U({\cal A}){\cal T} = [U({\cal W}^-) U({\cal Y}^-) ]{\cal T}\\
 &=& {\cal W}^- U({\cal W}^-)U({\cal Y}^-){\cal T} \oplus U({\cal Y}^-) {\cal T}.
\end{eqnarray*}
Since no nonzero element of ${\cal W}^- U({\cal W}^-)U({\cal Y}^-){\cal T}$ is a
lowest weight vector  
of ${\cal Y}$, ${\cal K}$ is the subspace of $\cal H$ of lowest
weight vectors of the abelian Lie algebra ${\cal Y}$.

In order to describe the lowest weight vectors explicitly
consider $S(\hat{\frak h}^-_{\Bbb Z})$ as a polynomial algebra on the
generators $x_i= w(-i)$, $z_i=r(-i), i >0$, so that
$S(\hat{\frak h}^-_{\Bbb Z}) = \Bbb C[x_i]_{i >0} \otimes \Bbb C[z_i]_{i >0}$.
The elements $w(i)$, $i \in \Bbb Z$ act on $S(\hat{\frak h}^-_{\Bbb Z})$ 
via multiplication, with $w(i)\cdot {\bf 1}= 0$ if $i>0$,
$[w(k), w(j)]=0$ and $[w(k), r(j)]= k \delta_{k+j,0} \la w, r \ra$. Thus
the element $w(k)$ acts on $\Bbb C[z_i]_{i >0}$ as the differential operator
$k\la r,w\ra \partial / \partial z_k$ for all $k >0$. 

By definition 
of the $K_i$, the lowest weight 
vectors of the ${\cal Y}$-module $\cal H$ are the lowest 
weight vectors of the above actions of the $w(i)$, $i \in \Bbb Z$, on 
$\Bbb C[x_i]_{i >0} \otimes \Bbb C[z_i]_{i >0}$. Since the action of the 
$w(i)$, $i \in \Bbb Z$,
commute with the elements $\Bbb C[x_i]_{i >0}$, the lowest weight vectors
are determined by the
elements $q  \in \Bbb C[z_i]_{i >0}$
satisfying $\partial / \partial z_k q = 0$ for all $k >0$, so
$q$ is a constant. Therefore the lowest weight vectors of the action of 
${\cal Y}$ on
$S(\hat{\frak h}^-_{\Bbb Z})$ correspond to the elements 
$\Bbb C[x_i]_{i >0} $. 
Thus 
${\cal K} = Ve^r \oplus [V \otimes (\Bbb C[x_i]_{i >0} \backslash \Bbb C) e^r] $.
Furthermore $V \otimes (\Bbb C[w(-i)]_{i >0} \backslash \Bbb C) e^r \subset
\mbox{rad}( \cdot, \cdot )_0$, and the form on $Ve^r$ is nondegenerate.
$\square$

Let ${\cal S}= {\cal W}^-U({\cal W}^-)U({\cal Y}^-)  {\cal T} \subset {\cal H}$.
This is called the space of ``spurious vectors'' in \cite{Tho}.
It follows from this definition (and the Poincare-Birkhoff-Witt theorem)
that $\cal H = S \oplus K$.

\begin{lemma} The associative algebra generated by the elements $L_i$
for $i > 0$ is generated by elements mapping ${\cal S}_{(1)}$ into $\cal S$.
\label{l:55}\end{lemma}

\noindent {\it Proof}\/: This is exactly the same argument as in
\cite{Tho}. First we will show that $L_1$ and $L_2 + {3 
\over 2} L_1^2$ have this property. Any $s \in {\cal S}$ can be written 
$$s = L_{-1} f_1 +L_{-2}f_2$$ where $ f_1, f_2 \in \cal H$ since any $
L_{m}, m<0$ can be written as a polynomial in $L_{-1}$ and $L_{-2}$. 
Furthermore $L_0s = s$ if and only if
$$L_0L_{-1}f_1 + L_0L_{-2}f_2
       = L_{-1}f_1 +L_{-2}f_2,
$$ and we may assume $ L_0 f_1 =0 , L_0f_2 = -f_2$.

Thus $s \in {\cal S}_{(1)}$ and $s =L_{-1}f_1 +L_{-2}f_2$ and $L_0f_2 
= -f_2 $ and $L_0 f_1 =0$.
Now we compute
\begin{eqnarray*}
L_1s &=& L_1L_{-1}f_1
+ L_1L_{-2}f_2\\
                 &= &L_{-1} L_1f_1 +2L_0 f_1 + 3L_{-1} f_2 + L_{-2}L_1 f_2 
\end{eqnarray*}
and this is in $\cal S$. Furthermore
\begin{eqnarray*}
(L_2 &+&{3 \over 2}L_1L_1) s  \\
       &=&L_2 L_{-1}f_1 + {3\over2}L_1L_1L_{-1}f_1 + L_2 L_{-2}f_2
        +{3\over2} L_1L_1L_{-2} f_2\\
       &=&L_{-1}L_2 f_1  + 3L_1f_1 + {3\over2}L_1L_{-1}L_1f_1 +L_{-2} L_2f_2\cr
       &+& 4L_0f_2 + {26\over12}6f_2 +{3\over2}L_1L_{-2}L_1 f_2 
        + {3^2\over2} L_1L_{-1}f_2\\
       &=&L_{-1}(L_2+{3\over2}L_1L_1)f_1 +L_{-2}(L_2 + {3\over2}L_1L_1)f_2
        + 9L_{-1}L_1f_2 .
\end{eqnarray*}
The above is a spurious vector (it contains $L_{-i}$ with $i >0$). Note
that $D= 26$ is necessary for this computation to work. Since $L_1$
and $L_2 +{3/2}L_1^2$ generate the algebra generated by the $L_i$, where
$i >0$, the lemma is proven.  $\square$

\begin{lemma} ${\cal P}_{(1)}$ is the direct sum of 
${\cal T}_{(1)}$ and  ${\cal N}_{(1)}$.
\label{l:56}\end{lemma}

\noindent{\it Proof}\/: Let $p \in {\cal
P}_{(1)}$. Then  $p =k + s$ where $k \in {{\cal K}}_{(1)}$ and $s \in
{\cal S}_{(1)}$; the decomposition is unique. By the preceding lemma
a generator $u$ (that is, $L_1$ or $L_2 + {3 \over 2} L_1^2$) of ${\cal W}^+$
satisfies $0=up = us+uk \in {\cal S} \oplus {\cal K}$. Thus $us=uk=0$,
and we see that
$s$ is annihilated by ${\cal W}^+$.
We conclude that $k \in {{\cal K} \cap {\cal P} = {\cal T}}$  
and $s \in {\cal S}_{(1)} \cap {\cal P}$. Since $\cal S$ is orthogonal
to $\cal P$,
$s$ must be an element in the radical of the form, $s \in {\cal
N}_{(1)}$. We conclude  
that ${\cal P}_{(1)} = {\cal T}_{(1)} \oplus {\cal N}_{(1)}$.$\square$

Theorem~\ref{thm:ng} now follows for $r \neq 0$ because 
Lemma A.4 and Lemma A.5 imply $Ve^r \approx {\cal T}$ so
$V_{(1-\la r,r \ra /2)} e^r 
\approx {\cal T}_{(1)}$, and Lemma~\ref{l:56} 
shows that ${\cal T}_{(1)}$ is
isomorphic to the quotient ${\cal P}_{(1)} /{\cal N}_{(1)}$. The
isomorphism is naturally a $G$-isomorphism. 

If $r =0$, first note that
$$( V \otimes V_{\cal L})_{(1)} = V_{(1)} 
\oplus (V_{(0)} \otimes (V_{\cal L})_{(1)}).$$
The subspace $(V_{\cal L})_{(1)}$ is two-dimensional, spanned 
by vectors of the form
$\alpha(-1)\iota(1), \beta(-1)\iota(1)$ where $\alpha, \beta$ span ${\cal L}$.
Furthermore, if $v \in  V_{(1)}$ then $L_nv =0$ for $n>1$ by the
condition that 
$V$ be $\Bbb N$-graded. Since $L_1 v \in V_{(0)}$, $L_1 v= c1$
for some constant $c$ and since $(L_1 v, 1)= (v, L_{-1}1)=0$ we have $c=0$.
It is easy to show that if $v \in V_{(0)} \otimes (V_{\cal L})_{(1)}$
then $L_nv =0$
for $n >0$. Thus the vectors in $ V_{(1)} \oplus (V_{(0)} \otimes
(V_{\cal L})_{(1)})$ 
are in ${\cal P}_{(1)}$, and the radical of the form restricted to this subspace 
is nondegenerate so Theorem~\ref{thm:ng} follows.
$\square$

\end{document}